\documentclass[12pt]{article}
\usepackage{geometry}
\usepackage{booktabs}

\usepackage{algorithm,algorithmic}

\usepackage{hyperref}

\usepackage{relsize}

\usepackage{bigints}

\usepackage{graphicx}

\usepackage{caption}

\usepackage[labelformat=simple]{subcaption}

\usepackage{multirow}

\usepackage{float}

\usepackage{array}

\usepackage{fullpage}

\usepackage{setspace}

\usepackage{paralist}

\usepackage{amsmath}
\usepackage{amssymb}
\usepackage{amsfonts}
\usepackage{amsthm}

\usepackage[left]{lineno}

\usepackage{tikz} \usetikzlibrary {plotmarks}

\newlength{\myrowheight}
\newcommand {\R}   {{\rm I\!R}}

\newcommand{\hf}{\frac12}

\newcommand{\E}{\vec{E}}

\newcommand{\s}{\vec{s}}

\newcommand{\B}{{\vec{B}}}

\newcommand{\curl}{\ensuremath{ \nabla \times\,}}

\newcommand{\bfB}{{\bf B}}

\newcommand{\bfE}{{\bf E}}

\newcommand{\bfe}{{\bf e}}

\newcommand{\bfb}{{\bf b}}

\newtheorem{definition}{Definition}

\usepackage{mathtools}
\usepackage{url}

\providecommand{\keywords}[1]{\textbf{\textit{Key words:  }} #1}
\providecommand{\msc}[1]{\textbf{\textit{AMS subject classification:  }} #1}

\begin{document}

\title{A Framework for the Upscaling of the Electrical Conductivity in the Quasi-static Maxwell's Equations}
\author{L. A. Caudillo-Mata*, E. Haber, L. J. Heagy, and C. Schwarzbach \\
lcaudill@eos.ubc.ca \\
Earth, Ocean and Atmospheric Sciences Department, \\ 
University of British Columbia, 4013-2207 Main Mall, \\
Vancouver, BC, Canada, V6T 1Z4}


{\let\newpage\relax\maketitle}

\begin{abstract}
Electromagnetic simulations of complex geologic settings are computationally expensive.  One reason for this is the fact that a fine mesh is required to accurately discretize the electrical conductivity model of a given setting.  This conductivity model may vary over several orders of magnitude and these variations can occur over a large range of length scales. Using a very fine mesh for the discretization of this setting leads to the necessity to solve a large system of equations that is often difficult to deal with.  To keep the simulations computationally tractable, coarse meshes are often employed for the discretization of the model. Such coarse meshes typically fail to capture the fine-scale variations in the conductivity model resulting in inaccuracies in the predicted data.  
In this work, we introduce a framework for constructing a coarse-mesh or upscaled conductivity model based on a prescribed fine-mesh model. Rather than using analytical expressions, we opt to pose upscaling as a parameter estimation problem.  By solving an optimization problem, we obtain a coarse-mesh conductivity model.  The optimization criterion can be tailored to the survey setting in order to produce coarse models that accurately reproduce the predicted data generated on the fine mesh.  This allows us to upscale arbitrary conductivity structures, as well as to better understand the meaning of the upscaled quantity. We use 1D and 3D examples to demonstrate that the proposed framework is able to emulate the behavior of the heterogeneity in the fine-mesh conductivity model, and to produce an accurate description of the desired predicted data obtained by using a coarse mesh in the simulation process.
\end{abstract}

\keywords{Numerical Homogenization, Maxwell's Equations, Simulation, Electrical Conductivity, Finite Volume, Geophysical Electromagnetic Methods}

\msc{Primary 78M40, 35K55; Secondary 65N08, 78A25, 86-08}
\section{Introduction}

Forward modeling of quasi-static Electromagnetic (EM) responses --- fields and fluxes --- significantly enhances our understanding of how variations in electrical conductivity impact EM data \cite{Ward1988,Weaver1999}. As a consequence, tools to simulate EM responses have been extensively used in a variety of scenarios, including medical applications \cite{Wang2014}, and geophysical applications such as mineral and hydrocarbon exploration, groundwater monitoring, and geotechnical and environmental investigations (cf. \cite{Oldenburg2007,Smith2013,Zhdanov2010}). 

Our particular interest is in the simulation of quasi-static EM responses over highly heterogeneous geologic settings, which are typically computationally expensive problems \cite{haberBook2014,Weaver1999}.  
Robust and accurate simulations of such settings require very fine meshes that are difficult, if not impossible, to work with as they translate into solving huge systems of equations (often on the order of tens of millions or even billions of unknowns).  
The simulation's cost is mainly due to the fact that the mesh must accurately capture the conductivity structures, which vary over a range of length scales and over several orders of magnitude. 
In terms of length-scales, the simulation domain may be on the order of tens of kilometers, topographic features may extend over several kilometers, the survey area may extend hundreds of meters, while geologic layers may vary on the scale of meters. 
The conductivity of geologic materials varies over more than ten orders of magnitude \cite{Telford1990}. This may result in areas with high conductivity contrasts in the geologic settings. All of these variations influence the behavior of the EM responses and must therefore be captured in the simulation; typically this is accomplished by using a very fine mesh.

Computational limitations often restrict EM simulations to a mesh that is coarser than the smallest length scale of significant features present in the geologic setting, therefore yielding inaccurate EM responses \cite{Avdeev2005, Borner2009}. 
Approaches to overcome the computational cost coarsen the mesh using adaptive, locally refined, or semi-structured meshes (cf. \cite{Haber2007a, Horesh2011,Key2011, Lipnikov2004, Schwarzbach2009}). These approaches have successfully produced accurate approximations to the EM responses at an affordable cost for some problems. However, these approaches face one major issue: the mesh must still capture the spatial distribution of the conductivity both inside and outside the region where we sample the EM responses to create data, thus restricting their ability to reduce the size of the system of equations to be solved.

Alternatively, upscaling or homogenization techniques seek to reduce the size of the system of equations to be solved by constructing  upscaled, homogenized or effective physical properties of the geologic medium, which are then used for simulation on a coarse mesh \cite{Wen1996}. These techniques have been extensively studied in the field of modeling flow in heterogeneous porous media, where they have been successfully used to drastically reduce the size of the system to be solved while producing accurate numerical solutions \cite{Durlofsky1998, Durlofsky2003,Farmer2002,MacLachlan2006, Moulton1998,Wen1996}.
Recognizing the success that upscaling techniques have had in fluid flow applications, we extend their use for application in EM  modeling.

Currently, the upscaling literature can be grouped in two main categories: analytical and numerical methods. 
Analytical methods use  averages to derive closed-form expressions for the upscaled partial differential equation (PDE) coefficients. They can be effective approaches when certain structural information of the medium to be homogenized is assumed.  The work in \cite{Wen1996, Farmer2002,Durlofsky2003,MacLachlan2004,Torquato2002,milton2002theory} provides an   overview of this type of analytical upscaling methods for application in simulating flow in porous media. 
Similar closed-form expressions can be obtained using effective medium theory, such as Archie's law which is used to homogenize the electrical conductivity of random heterogeneous materials (cf. \cite{BerrymanHoversten2013, KristenssonWellander2003, milton2002theory, Shafiro2000, Torquato2002}).  
These analytical and semi-analytical approaches have the advantage of being simple and affordable to compute; however, they are limited in the types of geologic media they can accurately describe. In contrast, numerical upscaling methods can handle more complex geologic media and they have proven to achieve accurate solutions comparable to those using traditional discretization methods, such as Finite Element (FE) or Finite Volume (FV), on the fine mesh. For applications  in the field of simulation fluid flow in porous media see \cite{Durlofsky1998, Durlofsky2003, Farmer2002, MacLachlan2006, Moulton1998} and references within. 

The purpose of this paper is to introduce a numerical framework to construct upscaled (coarse-mesh) electrical conductivity models of complex geologic settings for the quasi-static frequency domain EM problem. We tackle the problem of numerically constructing upscaled quantities in a manner that is fundamentally different than traditionally applied upscaling techniques, such as average-based techniques.  Those upscaling techniques do not consider the connection between the particular EM experiment configuration and the type of EM responses that the experiment produces, whereas this connection form the basis of the formulation for the upscaling problem we propose.  We formulate upscaling as a parameter estimation problem.  The proposed formulation has the following advantages.  First, it allows us to develop an application-specific framework for constructing coarse conductivity models.  Such coarse models aim to accurately reproduce the fine-scale EM responses by using a user-defined optimization criterion.  Second, it emphasizes the fact that there is no unique upscaled conductivity for different EM experiment configurations.  We will demonstrate this in Sections \ref{sec:upsWellLog} and \ref{sec:experiments}, where we show that using different frequencies in the experiment yield different upscaled conductivities, and similarly, different experiments result in different upscaled conductivities. 
Third, it can be tailored to construct either an isotropic or a fully anisotropic upscaled conductivity model depending on the context in which the upscaled quantity is to be used. That is, the framework is able to upscale arbitrary electrical conductivity structures.
This framework can be used to reduce the computational cost of EM forward modeling and to explore how local, small-scale conductivity variations affect the resulting EM responses.  To the best of our knowledge, the type of upscaling approach we propose has not been developed for Maxwell's equations within any other context. 

The proposed approach assumes a user-defined coarse mesh and a prescribed conductivity model discretized on a fine mesh, where the fine mesh sufficiently captures the significant conductivity variations present in the model.  A discrete parameter estimation problem is then formulated and solved on each coarse-mesh cell, yielding an upscaled conductivity quantity per coarse cell.  The aim is to construct upscaled quantities that emulate the effect of the fine-mesh conductivity inside each coarse cell on the EM responses of interest.  The upscaled conductivities are then assembled in a coarse-mesh electrical conductivity model that can be used for simulation.  We demonstrate the performance of the proposed upscaling framework using two examples and we compare it with traditional discretization techniques and other average-based upscaling procedures.   

The paper is organized as follows. In Section \ref{sec:upsCriteria}, we introduce the governing equations and formulate upscaling as a parameter estimation problem. In Section \ref{sec:upsWellLog}, we detail the upscaling method on a 1D example using well-log data.  In Section \ref{sec:ups3D}, a numerical upscaling framework in 3D is proposed. In Section \ref{sec:experiments}, numerical results for a representative synthetic 3D geologic case study over a mineral deposit are presented. Finally, in Section \ref{sec:conclusions} we discuss our findings and conclude the paper.

\section{An upscaling framework} \label{sec:upsCriteria}

This section introduces the mathematical framework that poses upscaling as a parameter estimation problem. 

\begin{figure}[bht!]
	\centering \includegraphics[width=.95\textwidth]{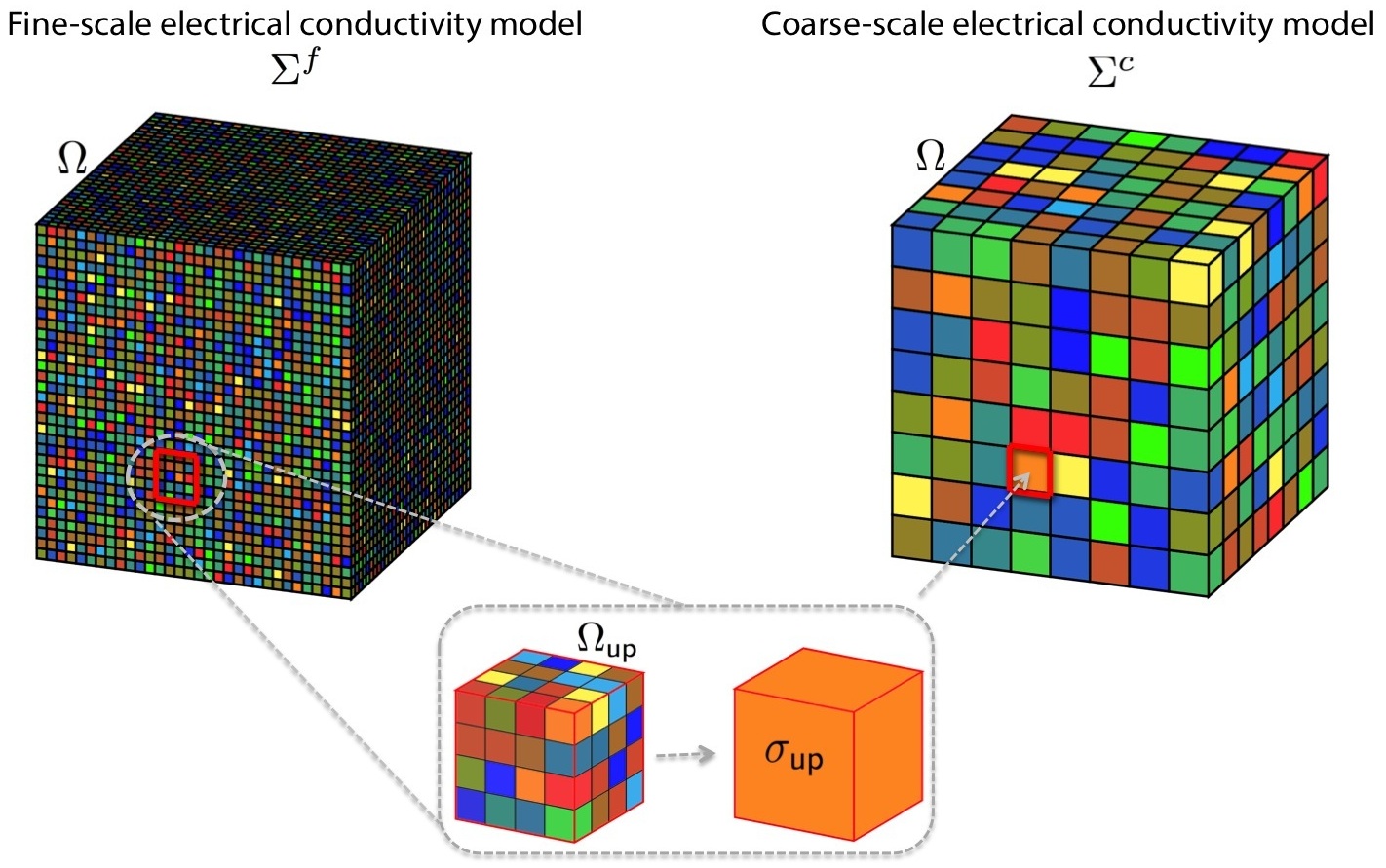}
	\caption{Upscaling process. Maxwell's equations are solved over the domain $\Omega$, where the electrical conductivity can vary over multiple length scales; the finest scale over which the conductivity varies defines the size of the mesh on which the model is discretized. The goal is to upscale the conductivity inside a subregion $\Omega_{\sf up}$ of $\Omega$ in order to construct a coarse-scale electrical conductivity model suitable for computation.  
	\label{fig:homog}}
\end{figure}

The EM problem is governed by Maxwell's equations.  Here, we consider the quasi-static Maxwell's equations in the frequency domain
\begin{equation}
\curl \E^f + i \omega \B^f = \vec{0}, \ \ \  \curl (\mu^{-1}\B^f) - \Sigma^f(\vec x) \E^f = \s;  \quad \vec x \in \Omega \subset \R^3,
\label{eq:EM_fine}
\end{equation}
where $\Omega$ is the domain, $\E^f $ is the electric field, $\B^f$ is the magnetic flux, $\s$ is the source term, $\omega$ is the angular frequency, $\mu$ is the magnetic permeability, and $\Sigma^f$ is a symmetric positive definite (SPD) tensor that represents the heterogeneous electrical conductivity.  The superscript $f$ refers to the fine scale, that is, the smallest scale (smallest spatial wavelength) over which the electrical conductivity varies significantly and that requires a fine mesh to be discretized on (Figure \ref{fig:homog}).

The PDE system \eqref{eq:EM_fine} is typically completed with the natural boundary conditions
\begin{equation}
	\mu^{-1}\B^f \times \vec{n} = \vec{0}, \label{eq:bcfine}
\end{equation}
where $\vec{n}$ denotes the unit outward-pointing normal vector to the boundary $\partial \Omega$ of $\Omega$. More general boundary conditions can be imposed as in \cite{Jin2002}; however, for the sake of simplicity, we limit the discussion to the boundary conditions \eqref{eq:bcfine}. 

Now, we are interested in solving numerically Maxwell's equations \eqref{eq:EM_fine}-\eqref{eq:bcfine} for highly heterogeneous geologic settings, where the heterogeneities may vary over multiple length scales. Since the conductivity is heterogeneous, possibly over a fine scale, a fine mesh is required for accurate simulations.  This leads to a large system of equations to be solved, making the simulation computationally expensive, or even intractable.  
One remedy to this problem is to replace the fine-scale conductivity structure, $\Sigma^{f},$ with a coarse-scale one, $\Sigma^{c}$. Given a $\Sigma^{c}$ that varies more slowly than $\Sigma^f$, Maxwell's equations can be discretized on a much coarser mesh, thus reducing the cost of the simulation. The process of obtaining $\Sigma^{c}$ from $\Sigma^{f}$ is referred to as {\em homogenization}
or {\em upscaling} (Figure \ref{fig:homog}).  

Note that we cannot simply replace the fine mesh with a coarse mesh in the simulation procedure. Defining a coarse mesh would necessarily require that several fine-mesh cells be captured in a single coarse-mesh cell (see central part of Figure \ref{fig:homog}). As a result, a single coarse cell may contain a heterogeneous fine-scale conductivity structure, meaning that the EM responses within this cell may be non-smooth.  Using a standard discretization technique (e.g. FE or FV) in this case will produce inaccurate approximations to the EM responses; such techniques assume a certain degree of smoothness in the function to be discretized.  To avoid this complication when replacing the fine mesh with a coarse mesh, we need to homogenize the fine-scale conductivity inside the coarse cells.  That is, for each coarse cell, we need to assign a representative quantity for the heterogeneous, fine-scale conductivity contained in it. Doing so, we ensure smoothness in the resulting EM responses inside each cell.

To derive a framework for the construction of $\Sigma^c$, we first consider the quasi-static Maxwell's equations on the coarse scale 
\begin{equation}
\curl \E^c + i \omega \B^c = \vec{0}, \ \ \ \curl (\mu^{-1}\B^c) - \Sigma^c \E^c = \s;  \label{eq:EM_coarse}
\end{equation}
where $\E^c, \B^c, \s, \omega, \mu$ are defined analogously as in \eqref{eq:EM_fine}. 
The superscript $c$ denotes dependency on the coarse scale, and the 
{\em coarse-scale conductivity model} is defined as follows 
\begin{eqnarray}
\label{condH}
	\Sigma^c(\vec x; \sigma_{\sf up}) = \left\{ \begin{matrix} \sigma_{\sf up}, & \mbox{if} \quad \vec x \in \Omega_{\sf up} \\
	\Sigma^f(\vec x), & \mbox{otherwise} 
	\end{matrix} \right.
\end{eqnarray}
where $\Omega_{\sf up} \subseteq \Omega$ is an {\em upscaling region} where we aim to homogenize the fine-scale conductivity (Figure \ref{fig:homog}), and $\sigma_{\sf up}$ is an {\em upscaled electrical conductivity} that aims to capture the effect of the fine-scale conductivity inside $\Omega_{\sf up}$ on the EM responses.  In this framework, the definition of $\sigma_{\sf up}$ depends on the context of a simulation. For example, $\sigma_{\sf up}$ may be given by a positive scalar, a real or even a complex tensor, depending on the complexity and aim of the simulation. We explore this idea further through examples in Sections \ref{sec:upsWellLog} and \ref{sec:experiments}.  

To simplify the presentation, we rewrite the fine-scale system \eqref{eq:EM_fine} and the coarse-scale system \eqref{eq:EM_coarse} as
\begin{equation}
 {\cal M}(\Sigma^f(\vec x)) \vec u^f = \vec q \quad {\rm and} \quad {\cal M}(\Sigma^c(\vec x; \sigma_{\sf up})) \vec u^c = \vec q,
\end{equation}
respectively.  Here, ${\cal M}$ represents the Maxwell operator, $\vec u^f = (\B^f,\E^f)^{\top}$ and $\vec u^c = (\B^c,\E^c)^{\top}$ represent the {\em fine-scale and coarse-scale EM responses}, respectively; and $\vec q$  encompasses the corresponding sources and boundary conditions \eqref{eq:bcfine}.  The fine-scale and coarse-scale EM responses are obtained by inverting the Maxwell operator, that is,
\begin{equation}
	\vec u^f(\Sigma^f(\vec x)) = {\cal M}^{-1}(\Sigma^f(\vec x))\vec q \quad {\rm and} 
	\quad  \vec u^c(\Sigma^c(\vec x; \sigma_{\sf up})) = {\cal M}^{-1}(\Sigma^c(\vec x; \sigma_{\sf up}))\vec q.
	\label{eq:fieldsFluxes}
\end{equation}

As previously stated, the constructed $\Sigma^c$ should be able to produce similar EM responses to the ones produced by $\Sigma^f$.  
We are typically interested in accurately simulating EM responses at certain locations within the simulation domain to produce predicted data at the receiver locations.
The simulated EM data at the receiver locations can be expressed as the action of a linear functional, ${\cal P}$, on the fine-scale and coarse-scale EM responses as follows
\begin{eqnarray}
\label{attributes}
d^f(\Sigma^f(\vec x)) = {\cal P} \vec u^f(\Sigma^f(\vec{x})) \quad {\rm and} \quad
d^c(\Sigma^c(\vec x; \sigma_{\sf up})) = {\cal P} \vec u^c(\Sigma^c(\vec{x};\sigma_{\sf up})).
\end{eqnarray}
Throughout this work, we refer to $d^f$ and $d^c$, as computed in \eqref{attributes}, as the {\em fine-scale and coarse-scale EM data}, respectively.  Note that when ${\cal P}$ equals the identity operator, \eqref{attributes} returns the EM responses in the entire simulation domain $\Omega$.

To conclude the construction of $\Sigma^c$, we need to define a criterion for choosing the `best' upscaled conductivity $\sigma_{\sf up}$ in the region $\Omega_{\sf up}$. That is,  we require a criterion 
able to construct a $\sigma_{\sf up}$ such that the coarse-scale data and the fine-scale data are similar.  We therefore propose the following definition:

\begin{definition} \label{def:homCond}
Let $\Sigma^f$ be the fine-scale electrical conductivity and let $\vec u^f$, given by \eqref{eq:fieldsFluxes}, be the resulting fine-scale EM responses (i.e., electric field and magnetic flux) for a given angular frequency $\omega$ and source (including boundary conditions) $\vec q$.
Let $d^f$ and $d^c,$ be some predicted fine and coarse-scale EM data as defined in \eqref{attributes}, then the upscaled electrical conductivity in the region $\Omega_{\sf up}$, denoted as $\sigma_{\sf up}^{*}$, is defined as the solution of the following parameter estimation problem
\begin{eqnarray}
\label{misfit}
\sigma_{\sf up}^{*} = {\rm arg}\min_{\sigma_{\sf up}} \, c(\sigma_{\sf up}) = \hf \big \|d^c(\Sigma^c(\vec x; \sigma_{\sf up})) - d^f(\Sigma^f(\vec x)) \big \|^{2}_2.
\end{eqnarray}
We refer to $c(\sigma_{\sf up})$ as the upscaling criterion.
\end{definition}
%
The proposed upscaling definition may look rather involved at first; however, an analogy can be drawn from the computation of an apparent conductivity in a Direct Current (DC) resistivity experiment. Although it would not be used for the purposes of simulation, the apparent conductivity can be considered as an upscaled quantity; for a given electrode geometry, the apparent conductivity is the homogeneous half-space conductivity that produces a response that is equivalent to the one observed. Within this upscaling context, the apparent conductivity corresponds to the quantity $\sigma^{*}_{\sf up}$, $\Omega_{\sf up}$ is the homogeneous earth (in this case, $\Omega_{\sf up}=\Omega$), $\vec{q}$ are the sources, $\vec{u}$ are the potentials, $\cal P$ projects the fields onto the receiving electrode locations, and the data $d^f$ and $d^c$ are the measured voltages. 

One aspect of definition \eqref{misfit} is that the upscaling region ($\Omega_{up}$) can be defined so that it takes into account the surrounding conductivity structure. That is, it corresponds to an `in situ' sampling of the electrical conductivity. This is demonstrated in the examples included in sections \ref{sec:upsWellLog} and \ref{sec:experiments}.

Another aspect of \eqref{misfit} is that the context of the simulation is considered in the construction of the upscaled quantity, both through the definition of source term and boundary conditions ($\vec{q}$), and choice of data of interest ($d^c$, $d^f$). For instance, the data of interest can be chosen among the electric or magnetic field or flux, or some combination of them. Additionally, the upscaled conductivity can be defined as isotropic or anisotropic, depending on the complexity of the setting and required accuracy of the coarse-scale simulation. Furtheremore, the upscaling criterion, $c(\sigma_{\sf up})$,  need not be based on least squares. The flexibility of the definition \eqref{misfit} allows it to be adapted depending on the expected complexity and intended use of the upscaled model.  All of these features provide a user-defined, application-specific framework.

Having an application-specific framework is important because it accounts for the fact that there is no unique upscaled conductivity suitable for all simulation contexts.
Indeed, constructing a different upscaling criterion by changing the data simulated or sources used to excite the system typically leads to a different upscaled conductivity.  The experiments we show in Sections \ref{sec:upsWellLog} and \ref{sec:experiments} demonstrate that by changing the frequency of the survey, it is possible to obtain different upscaled conductivity quantities. 
Although non-intuitive, different survey configurations have different sensitivity functions and sample the earth differently. Thus, the impact of heterogeneous conductivities may differ from experiment to experiment.  This same effect can be observed when considering DC resistivity surveys; for instance an apparent conductivity computed from a pole-dipole survey may be different than for a dipole-pole survey. 
We claim that definition \eqref{misfit} offers insights into the upscaling process. For example, it can explain why two different EM surveys that are conducted above the same area yield to different (upscaled) conductivities.

In the next sections we discuss the upscaling procedures in 1D and 3D.  


\section{Numerical upscaling in 1D} \label{sec:upsWellLog}

To demonstrate the use of the proposed upscaling framework, we begin with a 1D example. Our goal is to generate coarse-scale approximations of a well log electrical conductivity model for both a single frequency and a multi-frequency airborne loop-loop survey. Well log data have high resolution, with samples every few centimeters. However, if these conductivity models are to be used for earth-scale simulations or inversions, coarse-scale conductivity models, defined on the order of meters, are needed. In most practical applications, log data is simply averaged to obtain a homogenized conductivity, but as we see next, this may lead to serious errors. 

In this case study, we solve the parameter estimation problem defined in \eqref{misfit} numerically.  To do so, we use the discretize-then-optimize approach \cite{Haber2014}, thus we require both a careful discretization and an appropriate optimization method to solve the discrete version of \eqref{misfit}. 

\begin{figure}[!htb]
	\begin{subfigure}{0.5\textwidth}
		\centering
		\includegraphics[width=1\linewidth]{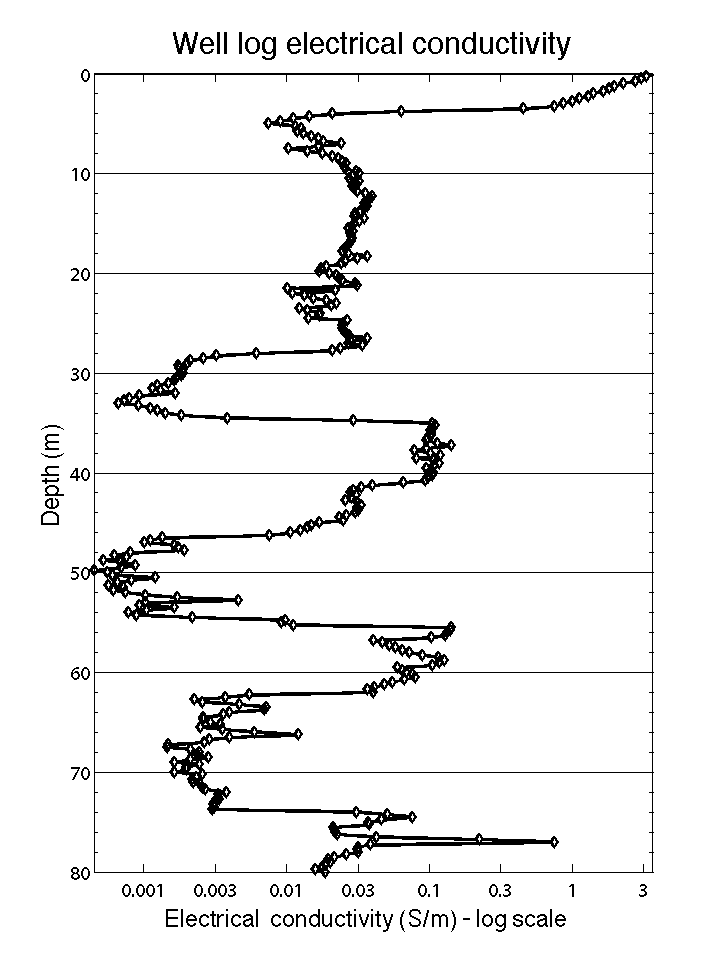}
		\caption{}
		\label{fig:fineC_unif_coarse}
	\end{subfigure}
	\begin{subfigure}{0.5\textwidth}
		\centering
		\includegraphics[width=1\linewidth]{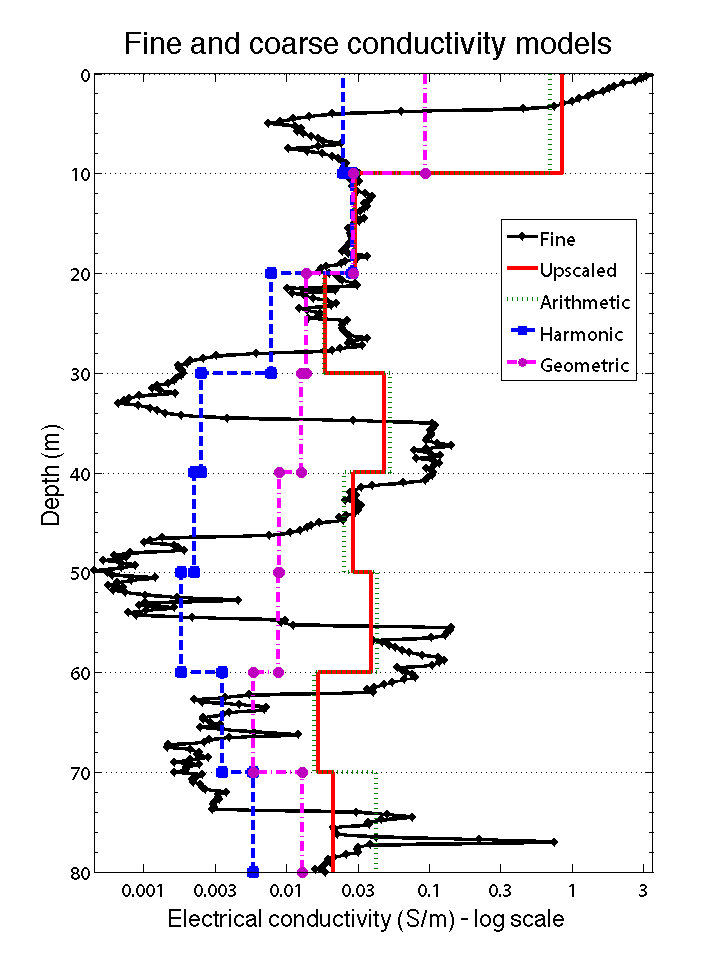}
		\caption{}
		\label{fig:singleF_upsC_unif}
	\end{subfigure}
    \caption{{Induction resistivity log: AA-05-01-096-11w4-0 from the McMurray/Wabiskaw Oil Sands deposit well log database. (a) Discrete fine-scale conductivity model. Each diamond represents a conductivity value on a uniform fine mesh of thickness 25 cm. Each straight line represents a coarse layer of a uniform mesh of thickness 10 m. The setup considers 320 fine layers and 8 coarse layers. (b) Resulting coarse-mesh conductivity models after applying four upscaling procedures: 1D numerical upscaling procedure (red solid line), and arithmetic (green dot line), geometric (magenta dot dash line), and harmonic (blue dash line) averages. }}
\end{figure}

We use a fine-scale conductivity model, $\Sigma^f$, given by an induction resistivity log from the
\href{http://www.ags.gov.ab.ca/publications/abstracts/SPE_006.html}{McMurray/Wabiskaw Oil Sands deposit}
well log public database \cite{wynne1994athabasca}. The McMurray formation is located in Northern Alberta, Canada, and the log used is shown in Figure \ref{fig:fineC_unif_coarse}.  
Observe that the interval over which  data were recorded is 80 m, thus we take this as the simulation domain, i.e., $\Omega = [0,80 ]$ m. In addition, observe that the electrical conductivity ranges over four orders of magnitude. The log chosen has 320 measurements total, with a measurement taken every 25 cm.  Hence, we defined a uniform fine mesh whose thickness is consistent with this scale.   

The first example uses an airborne survey configuration, with a frequency of 300 Hz and a horizontal coplanar arrangement for the source receiver-pair. The geometry of measurement system used
is taken from a standard Fugro survey. The source-receiver pair are located at a height of 40 m above the earth's surface and have a separation of 8.1 m. The source produces a magnetic field, which induces currents in the earth, producing secondary magnetic fields (H-field), which we measure at the receiver.

To construct a coarse-scale conductivity model that varies on the meter scale using the proposed upscaling framework, we need to choose: (a) a suitable coarse mesh, (b) the type of upscaled quantity to be constructed, and (c) an upscaling criterion.  

For the coarse mesh, we consider a uniform mesh nested in the fine mesh with 10 m thickness for each coarse layer. Hence, inside the simulation domain we have eight coarse layers, where each of them contains a large range of fine-scale conductivity variation to be upscaled, see Figure \ref{fig:fineC_unif_coarse}. The set of eight upscaled conductivities will form the discrete coarse-mesh conductivity model. 

We assumed the upscaled conductivity inside each coarse layer is given by a real positive scalar. In practice, this assumption is made when the resulting coarse-mesh conductivity model is to be used in a code that only handles isotropic variation of the conductivity.  A more general assumption is to consider the upscaled conductivity to be a tensor; however, this would require the forward modeling software to be capable of incorporating anisotropy, which is not always the case. We will elaborate further on 3D cases which incorporate anisotropy in Sections \ref{sec:ups3D} and \ref{sec:experiments}.

We are interested in simulating the behavior of the H-field at the receiver location, hence we take this as the EM data to be matched in the upscaling criterion \eqref{misfit}.  By doing so, we have specified the source term (including boundary conditions) and the data of interest; these define the necessary elements of the upscaling criterion for this example.  

We solve the discrete version of the parameter estimation problem \eqref{misfit} for each coarse layer separately. This yields the desired coarse-mesh conductivity model, which is shown with a red solid line in Figure \ref{fig:singleF_upsC_unif}. 

To compare our method to other 1D average-based upscaling methods, coarse-mesh conductivity models were constructed using arithmetic, geometric and harmonic averaging of the fine-mesh conductivity inside each coarse layer. The resulting coarse-mesh conductivity models are shown in Figure \ref{fig:singleF_upsC_unif}. 
Observe that the coarse-mesh conductivity model produced by the proposed upscaling procedure do not resemble an arithmetic, geometric or a harmonic average.  Rather than using the simple, context-independent averages (arithmetic, geometric or harmonic), the upscaling procedure accounts for the context of the simulation; it accounts for the source and surrounding conductivity structures in the construction of a homogenized quantity. Therefore, the proposed method gives an optimal \emph{in situ} prediction of the upscaled quantity, as defined in \eqref{misfit}, that is fundamentally different from the one given by the other average-based homogenization methods presented. 

To judge the quality of the various coarse-mesh conductivity models, we used each of them to forward model an H-field datum at the receiver location using the \href{http://gif.eos.ubc.ca/software/em1dfm}{EM1DFM} code on the coarse mesh. The EM1DFM code is based on the matrix propagation approach \cite{Farquharson2003}. Table \ref{tab:data1freqReg} shows the magnitude of the resulting H-field datum for each of the coarse-mesh models described, and the magnitude of the H-field datum obtained by forward modeling using the fine-mesh conductivity model.  The results in Table \ref{tab:data1freqReg} demonstrate that the proposed upscaling formulation constructed an optimal coarse-mesh conductivity model, in the sense of equation \eqref{misfit}, for the airborne survey configuration given, as it yields the smallest relative error in the approximation of the H-field datum of interest.  
The relative error is computed as the ratio of the absolute value of the difference in magnitude of the fine and coarse-mesh datum to the absolute value of the fine-mesh datum in magnitude. 
Note that by upscaling the conductivity model we reduced the fine-mesh problem size by 97\%.  

\begin{table}[!ht]
	\center
	\caption{{Magnitude of the magnetic field (H-field) datum and relative errors resulting from forward modeling using the fine-mesh and four coarse-mesh conductivity models. Note: $\%^{*}$ denotes percentage of primary field. }}\label{tab:data1freqReg}
	\begin{tabular}{  l  c  c  c  }
  		\toprule
  		\multicolumn{1}{p{3.0 cm}}{\centering Conductivity \\ model} & 
  		\multicolumn{1}{p{3.0cm}} {\centering Magnitude of \\ H-field ($\%^{*}$)} & 
  		\multicolumn{1}{p{3.0cm}} {\centering Relative error  \\ (percent)  } & 
  		\multicolumn{1}{p{3.0cm}} {\centering Number of \\ layers in mesh} \\        
  		\midrule  
  		Fine                             &    0.0549  & -----	 & 320   \\
  		Arithmetic                   &    0.0467  & 14.85 	 & 8 	\\
		Geometric 		&    0.0102  & 82.15 	 & 8 	\\
 		Harmonic 		&    0.0042  & 92.64 	 & 8 	\\
 		Upscaled		&    0.0535  &  6.29     & 8   	\\
  		\bottomrule
	\end{tabular}
\end{table}

The next example demonstrates the effect of considering multiple frequencies in the construction of the upscaled conductivity for the same survey configuration.

\subsection{Multi-frequency upscaling in 1D}  \label{multifreq}

For this example, we want to construct a coarse-mesh conductivity model to simulate magnetic field (H-field) measurements at five frequencies logarithmically equispaced in the range from $10$ Hz to $30,000$ Hz. In this section, we use the same airborne-style survey configuration; for the fine and coarse meshes, we use the same setup as described in the previous section.

\begin{figure}[!htb]
	\begin{subfigure}{0.5\textwidth}
		\centering
		\includegraphics[width=1\linewidth]{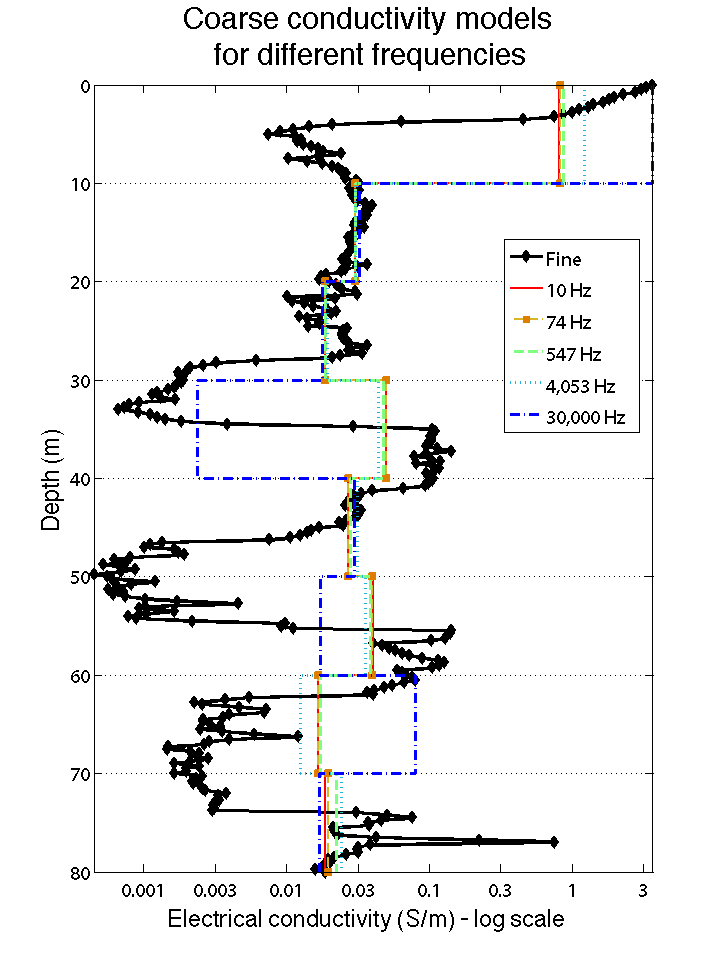}
		\caption{}
		\label{fig:multiF_upsC_unif}
	\end{subfigure}
	\begin{subfigure}{0.5\textwidth}
		\centering
		\includegraphics[width=1\linewidth]{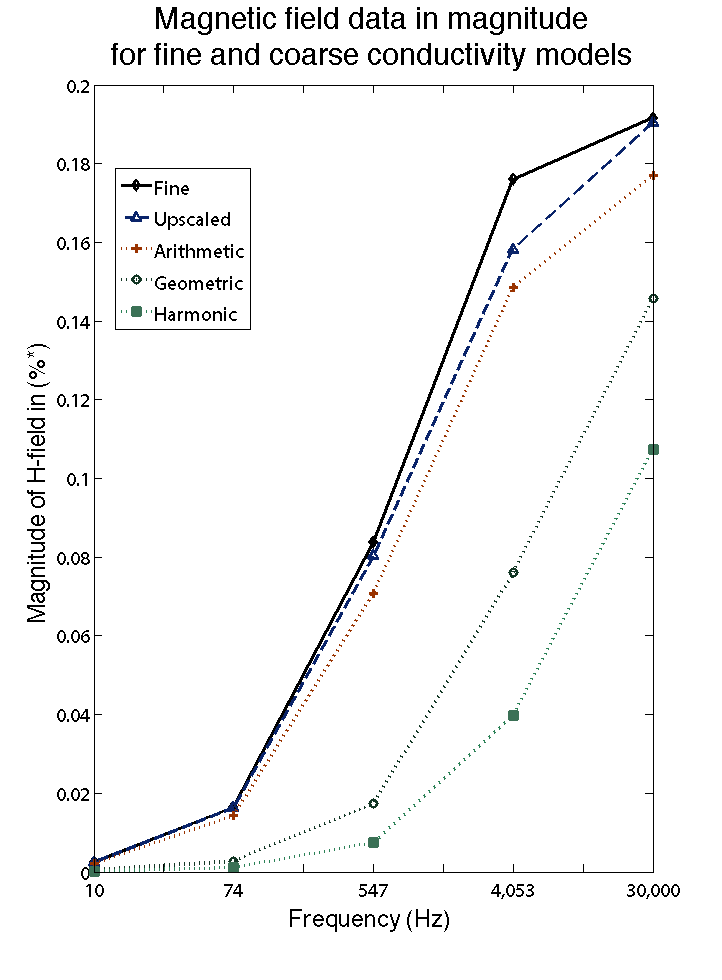}
		\caption{}
		\label{fig:multiF_mdata_unif}
	\end{subfigure}
	\caption{{(a) Coarse-mesh electrical conductivity models obtained by using the proposed upscaling method for different frequencies. The setup considers 320 fine layers and 8 coarse layers.  (b) Magnitude of magnetic field for each frequency, in \% of primary field ($\%^{*}$), resulting from forward modeling using the fine-mesh electrical conductivity model (black solid line), the different coarse-mesh conductivity models displayed in (a) (blue dash line), and the coarse-mesh models produced by using arithmetic (red plus dot line), geometric (gray circle dot line), and harmonic averages (green square dot line).}}
\end{figure}

Once again, we applied the proposed upscaling procedure by optimizing a discrete version of \eqref{misfit} for each of the eight coarse layers separately, and each individual frequency.  The resulting coarse-mesh conductivity models are shown in Figure   \ref{fig:multiF_upsC_unif}.  Observe that for each frequency we obtained a different coarse-mesh conductivity model. This follows from the fact that the upscaled conductivity model is tailored to match the magnetic field determined by the survey parameters.  These parameters influence the sensitivity of the magnetic field to the conductivity structure. Hence, varying any of these parameters alters how the conductivity structure is sampled. As a result, the upscaled conductivity model may take on different values depending on the experimental setting, demonstrating that the proposed upscaling approach provides a user-defined, application-specific framework.
These results also imply that the homogenized conductivity changes as a function
of frequency, demonstrating that frequency dependence on the coarse-scale may arise as a result of local, fine-scale heterogeneity. 

To conclude this example, we construct coarse-mesh conductivity models using arithmetic, geometric and harmonic averaging of the fine-mesh conductivity inside each coarse layer. For each of these coarse conductivity models, H-field data were then simulated on the coarse mesh using the given airborne survey configuration using the EM1DFM code (cf. \cite{Farquharson2003}). We compare the resulting H-field data with those computed using the fine-mesh conductivity model for each frequency in Figure \ref{fig:multiF_mdata_unif}.  The H-field data shown are given in percentage of the magnitude of the primary field. Table \ref{tab:multiF_rerrors} shows the relative errors obtained for each case.  
The relative error is computed as the ratio of the norm of the difference in magnitude of the fine and coarse-mesh data to the norm of the fine-mesh data in magnitude.
Observe that, once again, the proposed upscaling approach for constructing the coarse model produced better approximations to the desired H-field data than using the coarse models constructed by average-based upscaling procedures. 

\begin{table}[bht]
	\center
    \caption{{Relative errors for four coarse-mesh conductivity models and for five frequencies. }}\label{tab:multiF_rerrors}
	\begin{tabular}{  c  c  c  c  c  }
		\toprule
		\multicolumn{1}{c}{ } &
		\multicolumn{4}{c}{Relative errors for coarse-mesh conductivity models} \\ \cmidrule(r){2-5}
  		\multicolumn{1}{p{2.5 cm}}{\centering Frequency  \\ (Hz)} &
  		\multicolumn{1}{p{2.3 cm}}{\centering Arithmetic \\ (percent)} &
  		\multicolumn{1}{p{2.3 cm}}{\centering Geometric  \\ (percent)} &
  		\multicolumn{1}{p{2.3 cm}}{\centering Harmonic   \\ (percent)} &
  		\multicolumn{1}{p{2.3 cm}}{\centering Upscaled   \\ (percent)} \\
  		\midrule
  		10	    &     9.93 & 83.46	 & 93.34  &   0.64  \\
  		74	    &    12.76 & 83.36   & 93.23  &  2.78  \\
  		547     &    15.55 & 80.85 	 & 92.05  & 8.21   \\
		4,053	&    15.80 & 66.83 	 & 85.09  &11.76   \\
 		30,000	&     9.59 & 32.55 	 & 55.30  & 0.84   \\
  		\bottomrule
	\end{tabular}
\end{table}

Although 1D problems present minimal computational bottlenecks, if well log data are to be used to construct an earth-scale 3D conductivity model, or incorporated into an inversion, then using the proposed upscaling method, as shown in this section, provides a practical technique to reconcile these scales.

This case study illustrates the general principle behind the upscaling framework by showing its performance in a couple of examples using 1D well log data.  
However, applying the current upscaling formulation \eqref{misfit} to a general 3D setting is not practical. It requires simulating the data to be matched \eqref{attributes} on the fine mesh, which can be computationally demanding.  In the next section, we address the challenge of creating a practical upscaling approach for a 3D setting.


\section{Numerical upscaling framework in 3D} \label{sec:ups3D}

Applying the upscaling method used in the previous section for the 1D example is impractical to extend to the 3D case, as it requires a full simulation over the fine scale model, which is extremely expensive.  Therefore, in this section we propose a strategy for adapting the upscaling framework discussed in Section \ref{sec:upsCriteria} for practical application to 3D settings. 
To create a pragmatic method, we extended the methodology proposed by \cite{Durlofsky1998, Durlofsky2003} for the field of simulating flow in heterogeneous porous media to EM modeling. That is, we apply the upscaling procedure {\em locally}. In our case, this means that for each coarse-mesh cell, we locally solve a parameter estimation problem to construct an upscaled conductivity. Doing so cell by cell, potentially in parallel, yields the desired coarse-mesh conductivity model.
This approach enables us to solve several smaller problems rather than a single large one, making this procedure suitable for tackling large scale EM problems.

\begin{figure}[htb!]
 	\centering \includegraphics[width=0.95\textwidth]{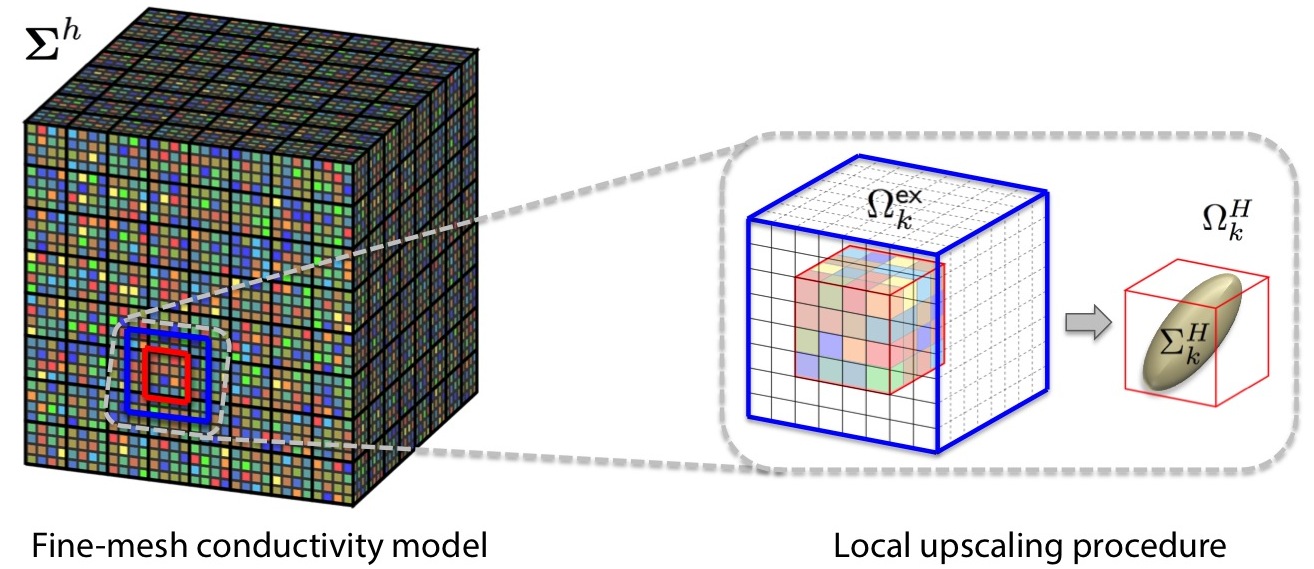}
 	\caption{{Local upscaling procedure for 3D settings.  Left: fine-mesh electrical conductivity model and example of nested meshes setup.  Right: extended domain ($\Omega_k^{\sf ex}$) for a given coarse-mesh cell ($\Omega^{H}_{k}$) and resulting anisotropic upscaled electrical conductivity (${\bf \Sigma}^H_k$). }
 	\label{fig:localUpscaling}}
\end{figure}

We now discuss, in detail, how to locally apply the upscaling framework. 
We assume that a given fine-scale conductivity model is discretized at the cell-centers of a 3D fine mesh, ${\cal S}^{h}$. The fine mesh sufficiently captures the significant conductivity variations in the model. We denote the discrete fine-mesh conductivity model as $\boldsymbol{\Sigma}^h$. 
We aim to construct a coarse-mesh conductivity model, $\boldsymbol{\Sigma}^{H}$, that is also discretized at the cell-centers of a user-chosen 3D coarse mesh, ${\cal S}^{H}$. Typically, ${\cal S}^{H}$ is much coarser than ${\cal S}^{h}$.
Throughout this section, the superscripts $h$ and $H$ denote dependency on the fine and coarse meshes, respectively.
The fine and coarse meshes are a union of $n$ fine cells and $N$ coarse cells, respectively. That is, ${\cal S}^{h} = \cup_{i=1}^n \Omega^{h}_{i}$ and ${\cal S}^{H} = \cup_{k=1}^N \Omega^{H}_{k}$, where $N \ll n$. 
For simplicity, we also assume that the meshes are nested, that is ${\cal S}^{H} \subseteq {\cal S}^{h}$; however, the argument presented here can be extended to include more general mesh setups. A sketch of the setup described is shown on the left-hand side of Figure \ref{fig:localUpscaling}.

Since the goal is to apply the upscaling procedure on each coarse-mesh cell independently (locally), we need to identify: (a) the upscaling region, (b) the type of upscaling quantity to be constructed, and (c) the data to be matched in a local version of the parameter estimation problem \eqref{misfit}. We discuss each of these for a single coarse cell below.

Consider a single coarse-mesh cell, $\Omega^{H}_{k}$. The upscaling region corresponds to $\Omega^{H}_{k}$, which is composed of the fine cells and the fine conductivity structure it encloses.  To construct an upscaled conductivity in $\Omega^{H}_{k}$ that takes into account the surrounding conductivity structure (i.e.,  to preserve the `in situ' behavior of the fine-scale conductivity), we embedded $\Omega_k^H$ in an extended domain, $\Omega_k^{\sf ex}$. This extended domain $\Omega_k^{\sf ex}$ includes $\Omega_k^H$ and a neighborhood of fine cells (and their corresponding conductivity values) around $\Omega_k^H$, see Figure \ref{fig:localUpscaling}.   The size of the extended domain is user-chosen.  We explore the effect of different extensions in the next section.

To better represent most of the existing heterogeneity inside the cell $\Omega_k^H$, we assume the upscaled conductivity to be constructed for this cell, ${\bf \Sigma}^H_k$, is fully anisotropic.  That is, ${\bf \Sigma}^H_k$ is a SPD matrix, which can be parametrized as 
\begin{equation}
\label{eq:sigmaParametrized}
{{\bf \Sigma}^H_k}(\sigma^k_1,\sigma^k_2,\sigma^k_3,\sigma^k_4,\sigma^k_5,\sigma^k_6) = \begin{bmatrix}
                  \sigma^k_{1} & \sigma^k_{4} & \sigma^k_{5} \\
                  \sigma^k_{4} & \sigma^k_{2} & \sigma^k_{6} \\
                  \sigma^k_{5} & \sigma^k_{6} & \sigma^k_{3}
                  \end{bmatrix}; \ \ \ \{\sigma^k_j\}_{j=1}^6 \in \R.
\end{equation}

According to Definition \ref{def:homCond} (Section \ref{sec:upsCriteria}), in order to construct ${\boldsymbol \Sigma}^{H}_k$ by solving a parameter estimation problem, we require some fine and coarse-scale data to be matched in \eqref{misfit}.  To generate such data, the Maxwell system should be excited by either a source or some boundary conditions, see \eqref{attributes}.
Since we apply the upscaling procedure locally, we assume that sources do not reside inside $\Omega_k^{\sf ex}$. Therefore, rather than choosing some local source(s) to induce EM responses, we assume that the system is excited by some non-homogeneous boundary conditions.  Such boundary conditions should reflect the behavior of the EM responses of interest across the boundary of $\Omega_k^{H}$, denoted as $\partial \Omega_k^{H}$.
In principle, the correct boundary conditions can be obtained numerically by solving the fine-mesh problem; however, they are impractical to compute.  One remedy for this problem, suggested in \cite{Durlofsky2003,Haber2014c}, is to use a set of linearly independent boundary conditions.  Note that using linear boundary conditions in the context of the EM problem can be appropriate to model the action of distant sources through $\partial \Omega_k^H$ as such action can be perceived as a `plane wave'.

\begin{figure}[htb!]
 	\centering \includegraphics[width=1.04\textwidth]{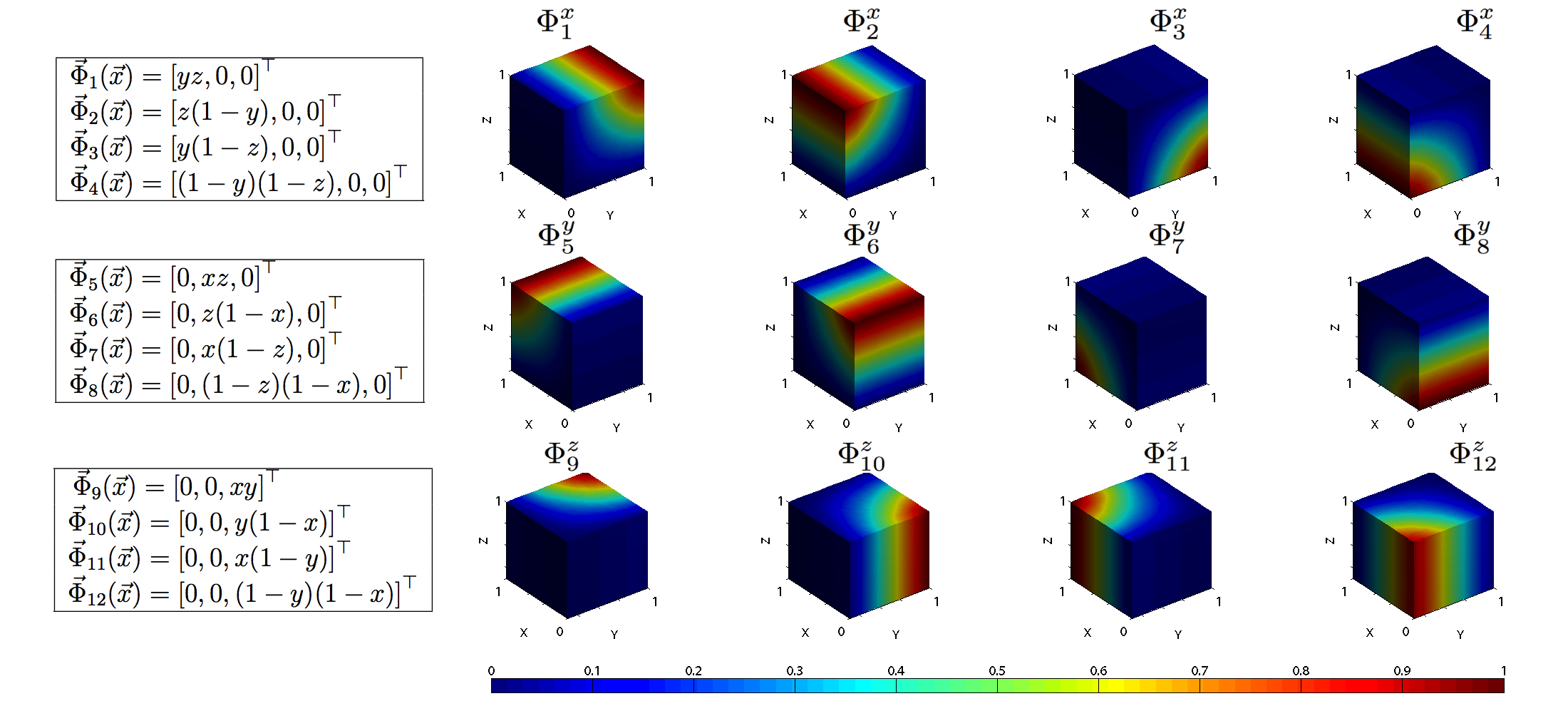}
 	\caption{{Left: analytical expressions for the set of linearly independent boundary conditions used to locally generate data in $\Omega_k^{\sf ex}$ (one per edge).  Right: plot of non-zero component functions for each $\vec{\Phi}_l$ on $\partial \Omega_k^{\sf ex}$.} }
 	\label{fig:linDecBC}
\end{figure}

We generate fine and coarse-scale data by locally exciting the Maxwell system  using a set of linearly independent boundary conditions, one per edge of $\Omega_k^{\sf ex}$. This yields the following set of twelve local problems 
\begin{eqnarray}
	\curl \E_l^k + i \omega \B_l^k &=& \vec{0}, \ \ \forall \vec{x} \in \Omega_k^{\sf ex}; \label{EM1_k} \\
	\curl (\mu^{-1} \B_l^k) - \Sigma^k(\vec x) \E_l^k &=& \vec{0}, \ \  \forall \vec{x} \in \Omega_k^{\sf ex}; \label{EM2_k} \\
	\E_l^k \times \vec{n} &=& \vec{\Phi}_l \times \vec{n}, \ \ \forall \vec{x} \in \partial \Omega_k^{\sf ex}, \ \ l=1,...,12, \label{BC_k}
\end{eqnarray}
where $\partial \Omega_k^{\sf ex}$ denotes the boundary of $\Omega_k^{\sf ex}$, and each $\vec{\Phi}_l$ is a vector function that takes the value 1 along the tangential direction to the $l$-th edge and decays linearly to 0 in the normal directions to the $l$-th edge. See Figure \ref{fig:linDecBC}.  This set of boundary conditions form the natural basis functions for edge degrees of freedom (cf. \cite{Haber2014c,Monk2003}), therefore they can be used to model general linearly varying EM responses.  Similar choices were proposed in \cite{Durlofsky2003,Efendiev2009} for the problem of simulating fluid flow in porous media, where the PDE model is the Poisson equation.
Different studies for flow in porous media have shown that different choices of boundary conditions lead to different upscaled quantities (cf. \cite{Durlofsky1998,Durlofsky2003}). 

To compute the fine-scale data for the $k$-th local problem, we use the fine-mesh conductivity contained in $\Omega_k^{\sf ex}$, denoted as ${\bf \Sigma}^h_k$. To compute the coarse-scale data for the $k$-th local problem, we assign ${\bf \Sigma}^H_k$ as the conductivity value of every fine-mesh cell inside $\Omega_k^{\sf ex}$.

Using the fine-mesh conductivity inside $\Omega_k^{\sf ex}$ (Figure \ref{fig:localUpscaling}), we can now forward model each of the twelve problems \eqref{EM1_k}-\eqref{BC_k} to obtain a set of EM responses.  To do so, we discretize each of these problems using traditional edge-based discretization methods, such as FE (cf. \cite{Jin2002}) or Mimetic Finite Volume (MFV) (cf. \cite{Hyman1998,Hyman1999,Hyman1999a,Haber2014}).  Edge-based discretization methods use staggered meshes that discretize $\E$ on the edges, $\B$ on the faces, and the PDE coefficients ($\Sigma$ and $\mu$) at the cell-centers of a given cell. 
In Section \ref{sec:experiments}, we show results where the discretization was done using the MFV method; however, an edge-based FE discretization method can be used as well.
After the discretization, we obtain a set of twelve discrete electric fields, $\bfE_k = \{\bfe^k_{1},\ldots,\bfe^k_{12}\}$, and a set of twelve discrete magnetic fluxes, $\bfB_k = \{\bfb^k_{1},\ldots,\bfb^k_{12}\}$.  

Next, we need to choose which fine and coarse data are to be matched and how the local version of the upscaling criterion~\eqref{misfit} will be formulated. We note that the heterogeneous fine-scale conductivity structure inside $\Omega_k^H$ generates non-smooth and potentially discontinuous responses on $\partial \Omega_k^H$. Since the role of the upscaled conductivity is to emulate the effect of the fine-scale conductivity inside $\Omega_k^H$ on the EM responses at $\partial \Omega_k^H$, it makes sense to consider the data to be either the integral of the electric field over the twelve edges of $\partial \Omega_k^H$, or the integral of the magnetic flux over the six faces of $\partial \Omega_k^H$.  We refer to these data as the total electric fields or total magnetic fluxes, respectively.
To be more specific, we define the {\em total electric field data} as
\begin{subequations}
\label{data3d}
\begin{eqnarray}
 d_{lm} = \int_{\Gamma_{m}} \E_{l}^k \cdot \vec{\tau}~d\ell; \quad l=1,\ldots,12, \ m=1,\ldots,12,
 \end{eqnarray}
and the {\em total magnetic flux data} as
\begin{eqnarray}
 d_{lj} = \int_{F_{j}} \B_{l}^k \cdot {\vec n}~dS;   \quad l=1,\ldots,12, \ j=1,\ldots,6,
 \label{eq:totalB}
 \end{eqnarray}
 \end{subequations}
where $\Gamma_{m}$ represents the $m$-th edge of $\partial \Omega_k^H$, $\vec{\tau}$ denotes the unit tangent vector to $\Gamma_{m}$, $F_{j}$ represents the $j$-th face of $\partial \Omega_k^H$, and $\vec{n}$ represents the unit normal vector to $F_{j}$. The data can be computed by numerically integrating \eqref{data3d} using the set of discrete fields $\bfE_k$ or fluxes $\bfB_k$, respectively.  Remember that $\bfE_k$ and $\bfB_k$ are discretized on the fine-mesh edges and fine-mesh faces inside $\Omega_k^{\sf ex}$, respectively.
As before, the choice of the data depends on the context of a given simulation, which we will demonstrate for a 3D example in the next section. The choice of the data to be matched during the local upscaling procedure significantly influences the construction of the upscaled conductivity; we have shown this in \cite{caudillo2014}.									

Finally, we formulate the local version of the discrete parameter estimation problem \eqref{misfit} to be solved in order to construct ${\bf \Sigma}^H_k$ as follows
\begin{eqnarray}
\begin{aligned}
{\bf \Sigma}^{H}_k = \mbox{} & \underset{ \hat{\bf \Sigma}^H_k  \in \R^{3 \times 3} }{ \text{arg min} } \
\frac{1}{2} \sum_{l=1}^{12} \left \|{\bf d}_l\big(\hat{\bf \Sigma}^H_k\big) - {\bf d}_l\big({\bf \Sigma}^h_k\big) \right \|_2^2 \\
& \text{subject to } \hat{\bf \Sigma}^H_k \text{ is SPD.}  \label{eq:optP}
\end{aligned}
\end{eqnarray}
Here, ${\bf d}_l\big({\bf \Sigma}^h_k)$ denotes the vector whose entries are the $l$-th total electric field data or total magnetic flux data \eqref{data3d} computed for the fine-mesh conductivity in $\Omega_k^{\sf ex}$ (i.e., ${\bf \Sigma}^h_k$); ${\bf d}_l(\hat{\bf \Sigma}^H_k)$ denotes the vector whose entries are the $l$-th total electric field data or total magnetic flux data \eqref{data3d} computed for the upscaled conductivity $\hat{\bf \Sigma}^H_k$ extended to $\Omega_k^{\sf ex}$. 
To solve the discrete optimization problem \eqref{eq:optP}, we first compute the gradients and sensitivities as discussed in \cite{Haber2014}, and then use the Projected Gauss-Newton method (cf. \cite{Kelley1999,Lin1999,Nocedal2006}). 

In summary, given a discrete fine-mesh conductivity model ${\boldsymbol \Sigma}^{h}$, we can compute an upscaled conductivity (${\boldsymbol \Sigma}^{H}_k$) in a coarse cell ($\Omega_k^H$) by the following steps:
\begin{enumerate}
\item Choose the size of the extended local domain ($\Omega_k^{\sf ex}$) where the coarse cell ($\Omega_{k}^H$) is embedded. See Figure \ref{fig:localUpscaling}.
\item Forward model the twelve problems \eqref{EM1_k}-\eqref{BC_k} defined on $\Omega_k^{\sf ex}$ to obtain the sets of discrete electric fields ($\bfE_k$) and magnetic fluxes ($\bfB_k$).  To do so, we use a traditional edge-based discretization method, such as FE or MFV. The discretization is done using the fine mesh inside $\Omega_k^{\sf ex}$ and its corresponding fine-mesh conductivity information ${\boldsymbol \Sigma}^{h}_k$.  
\item Choose the data set to be matched, that is, total electric field data or total magnetic flux data as defined in \eqref{data3d}, according to the context of the given simulation. Compute the data ${\bf d}_l\big({\bf \Sigma}^h_k\big)$ using the discrete fields and fluxes obtained in the previous step.
\item Optimize the constrained parameter estimation problem \eqref{eq:optP} to obtain the desired anisotropic upscaled conductivity ${\bf \Sigma}^{H}_k$.  Note that solving such an optimization problem involves performing steps 2 and 3 (for the upscaled conductivity model $\hat{\bf \Sigma}^H_k$) to compute ${\bf d}_l\big(\hat{\bf \Sigma}^H_k\big)$ at each iteration.
\end{enumerate}

Observe that all the calculations are done on $\Omega_k^{\sf ex}$, where each optimization problem is small and can be solved quickly. Furthermore, since the problem defined for each coarse cell is independent, the upscaling procedure can be done in parallel.
Once all the upscaled conductivities are computed, we assemble the coarse-mesh conductivity model.  Using this coarse conductivity model, the original problem can then be discretized and solved on the coarse mesh.  In particular, we have obtained a system of equations that is substantially smaller than the original one.

 
\section{Numerical upscaling results in 3D} \label{sec:experiments}

In this section, we demonstrate the performance of the proposed numerical upscaling framework, introduced in Section \ref{sec:ups3D}, for a 3D example.  Our goal is to construct anisotropic, coarse-mesh electrical conductivity models which can be used to simulate magnetic data for a large-loop EM survey over a synthetic model of a mineral deposit.  The results are compared to the ones obtained on a fine mesh. We also compare our results with some other simple averaging approaches.

\begin{figure}[ht]
  \centering
  \begin{tabular}[c]{cc}
    \begin{subfigure}[c]{0.48\textwidth}
      \includegraphics[width=\textwidth]{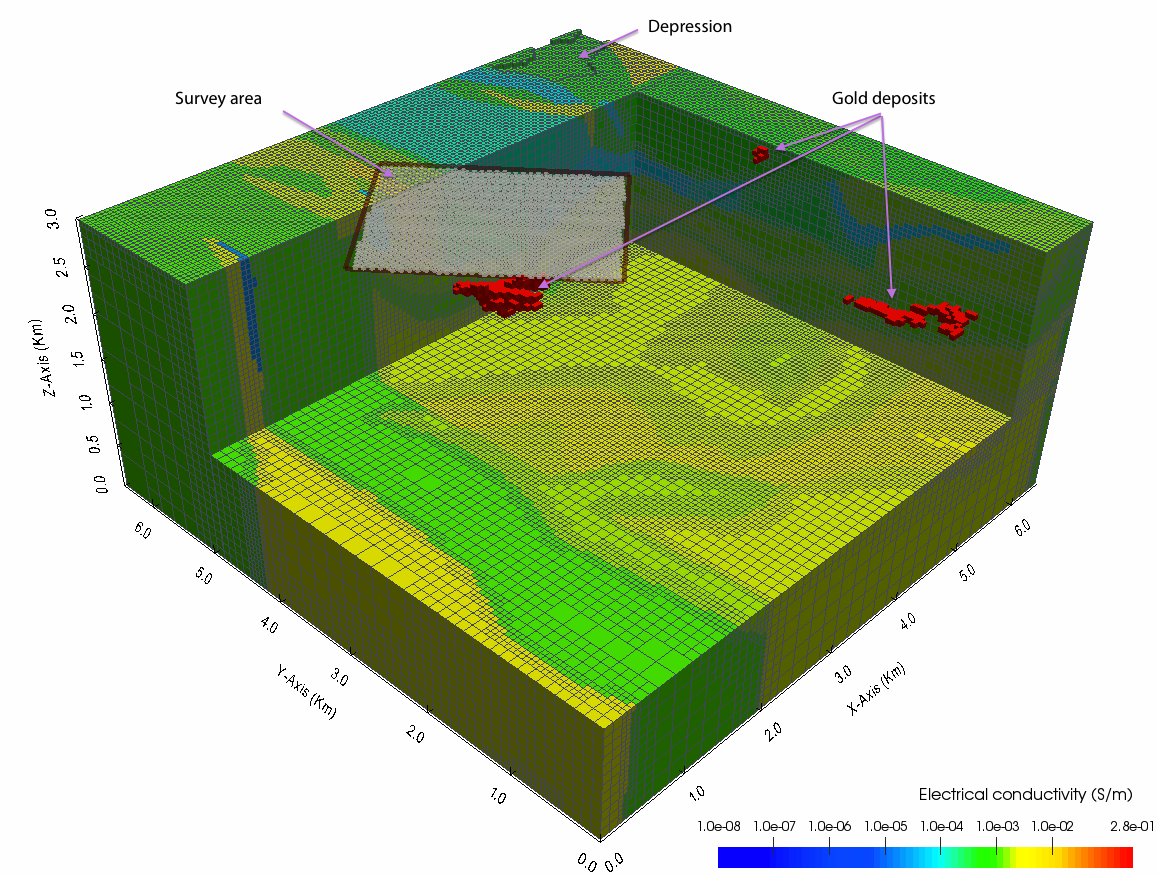}
      \caption{}
      \label{fig:referenceModel}
    \end{subfigure}&
    \begin{subfigure}[c]{0.48\textwidth}
      \includegraphics[width=\textwidth]{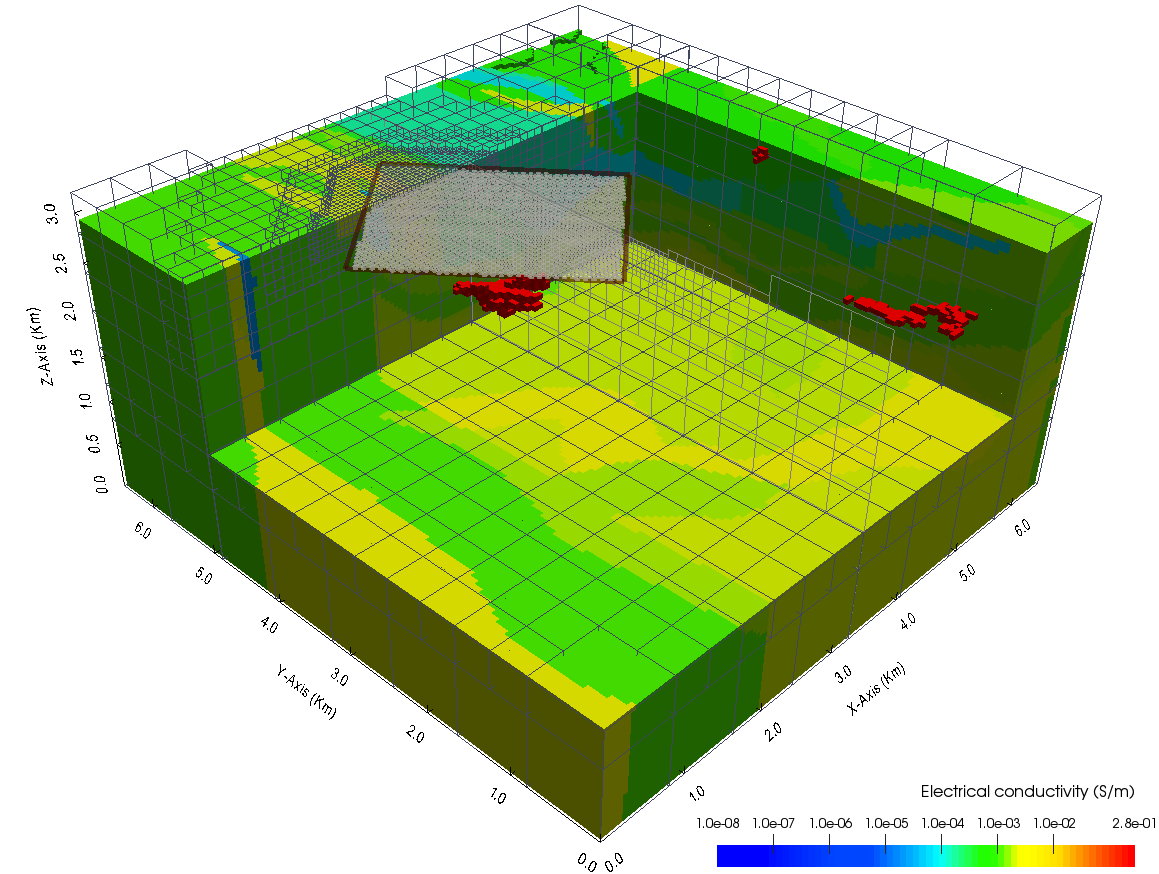}
      \caption{}
      \label{fig:CoarseLalorSetting}
    \end{subfigure}\\
  \end{tabular}    
  	\caption{Subsurface part of the synthetic electrical conductivity model on a fine and on a coarse OcTree mesh, and large-loop EM survey setting. (a) Model discretized on a fine OcTree mesh (546,295 cells). The conductivity varies over six orders of magnitude. (b) Model with an overlaying coarse OcTree mesh (60,656 cells). The coarse OcTree mesh maintains the same cell size as the fine mesh in the survey area and gradually increases the cell size for the rest of the domain. We use the proposed 3D upscaling framework to construct an anisotropic conductivity model for the coarse mesh.}
  	\label{fig:condModels}
\end{figure}

For this example, we constructed a synthetic conductivity model based on an inversion model obtained by \cite{Yang2014} using field measurements over the Canadian Lalor mine. 
The Lalor mine targets a large zinc - gold - copper deposit that has been the subject of several EM surveys (see \cite{Yang2014} and references within).  In particular, this conductivity model contains three highly conductive units which represent the deposit.  The model has  non-flat topography that requires an adaptive mesh to be faithfully represented,  see Figure \ref{fig:referenceModel}.  The domain of the model extends from 0 to 6.5~km along the $X$, $Y$ and $Z$ directions, respectively. The model contains air and the subsurface. We assumed a constant conductivity value of $10^{-8}$~S/m in the air. The subsurface is composed of 35 geologic units, including the three conductive units comprising the deposit (Figure \ref{fig:referenceModel}). The subsurface conductivity values range from $1.96 \times 10^{-5}$~S/m to $0.28$~S/m. 
We discretized the conductivity model at the cell centers of a fine OcTree mesh.
The OcTree mesh used has cell sizes of $(50~m)^{3}$  within the area where the EM survey will be conducted and at the interfaces of the features in the model where the conductivity varies. The rest of the domain was padded out with coarser OcTree cells. The fine OcTree mesh is shown in Figure \ref{fig:referenceModel}; it has 546,295 cells.  More details about OcTree meshes can be found in \cite{Haber2007a,Horesh2011}. 

We considered a large-loop EM survey for this example. A rectangular transmitter loop, with dimensions 2~km $\times$ 3~km, was used.  The transmitter is placed on the Earth's surface above the largest deposit, as shown in Figure \ref{fig:referenceModel}, and operates at the frequencies of 1 and 20~Hz.  Inside the transmitter loop, we placed a uniform grid of receivers that measure the three components of the magnetic flux.  The measurement locations are separated by 50~m along the $X$ and $Y$ directions, respectively.  We refer to the rectangular area where the source and receivers are located as the {\em survey area}. In this example, we are interested in accurately simulating magnetic flux data within this area. To get an estimate for the proper cell size, we considered the largest background conductivity value (4.5$\times 10^{-3}$~S/m) and calculated skin depths of 7,461 and 1,668~m for the frequencies of 1 and 20~Hz, respectively.  Thus, using cells of size 50, 100, and 200~m should be sufficient to capture the decaying nature of the magnetic flux on the finer OcTree mesh. Observe from Figure \ref{fig:referenceModel} that the survey area was embedded into a much larger computational  domain to reduce  the effect of the imposed boundary conditions \eqref{eq:bcfine}, which replace the true decay of the fields towards infinity.

In order to construct an anisotropic coarse conductivity model using the upscaling methodology introduced in Section \ref{sec:ups3D}, we need to choose the following parameters:
(a) a suitable coarse mesh,  (b) the size of the extended local domains, and (c) the data to be matched in the upscaling criterion.

As a coarse mesh, we considered an OcTree mesh that is nested within the fine mesh. 
Since we are interested in accurately simulating magnetic flux data in the survey area, the OcTree mesh was designed to maintain the fine mesh resolution $(50~m)^3$ inside the area of interest, whereas the rest of the domain was filled with increasingly coarser cells. In total it contains only 60,656 cells, that is, roughly $10\%$ of the number of cells in the fine OcTree mesh.  Figure \ref{fig:CoarseLalorSetting} shows the coarse OcTree mesh. 
Observe from Figure \ref{fig:CoarseLalorSetting} that the coarse OcTree mesh was not refined outside the survey area where a large conductivity contrast is present in the model. For example, at the interface between the highly conductive gold units and the more resistive background, and at the air-Earth interface which is not flat. These interfaces are not represented in the coarse OcTree mesh. We challenge the upscaling procedure with the large conductivity difference across the interfaces. 

Next, we need to choose the size of the extended domain to solve the local problems, that is, the number of fine-mesh padding cells by which we extend every coarse cell to be upscaled as shown in Figure \ref{fig:localUpscaling}.  To investigate the effect of this size on the resulting upscaled conductivity, we performed the extension using two sizes: 4 and 8 padding cells.

Since the receivers at the survey area measure magnetic flux data, we considered the total magnetic flux \eqref{eq:totalB} through each of the faces of the coarse cell as the data to be matched by the upscaling criterion. By doing so, we have connected the upscaling criterion to the large-loop EM survey.

We applied the upscaling framework introduced in Section \ref{sec:ups3D} by solving the constrained parameter estimation problem for each of the coarse cells separately, each individual size of the extended domain, and each individual frequency.  This yields two anisotropic, coarse-mesh conductivity models.

\begin{figure}[ht]
  \centering
  \begin{tabular}[c]{ccc}
   \toprule
    \multicolumn{1}{c}{} &
	\multicolumn{2}{c}{Magnitude of the magnetic flux for the large-loop EM survey setting} \\ \cmidrule(r){2-3}
	&
	\multicolumn{1}{p{6 cm}}{\centering Reference solution }   & 
	\multicolumn{1}{p{7 cm}}{\centering Upscaled solution (8 padding cells)}  \\
   \midrule
   \rotatebox{90}{1 Hz} &
    \begin{subfigure}[c]{0.45\textwidth}
      \includegraphics[width=\textwidth]{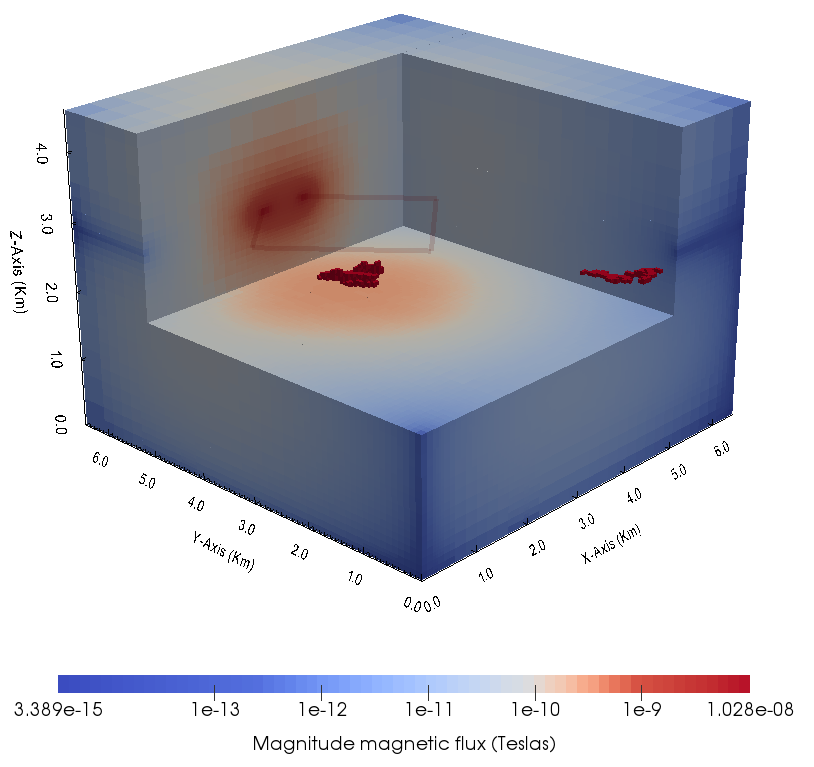}
    \end{subfigure}&
    \begin{subfigure}[c]{0.45\textwidth}
      \includegraphics[width=\textwidth]{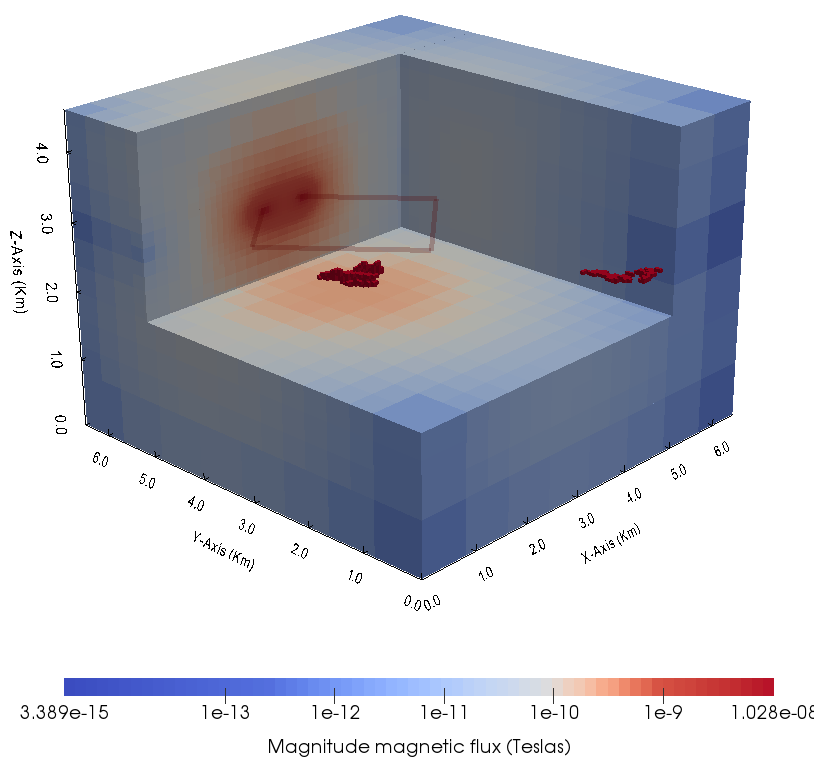}
    \end{subfigure}\\
    \rotatebox{90}{20 Hz} &
    \begin{subfigure}[c]{0.45\textwidth}
      \includegraphics[width=\textwidth]{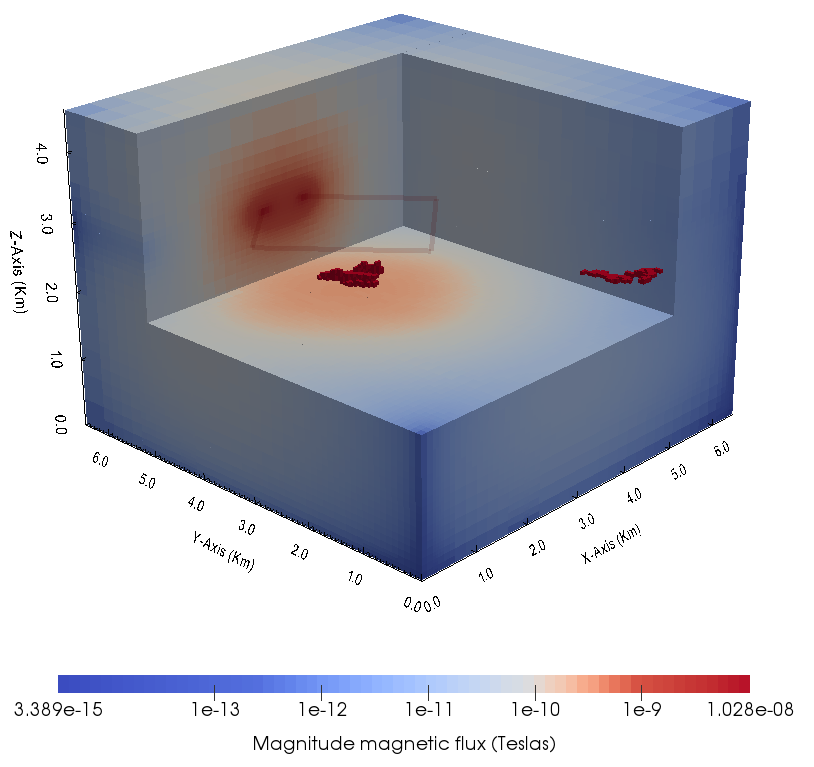}
    \end{subfigure}&
    \begin{subfigure}[c]{0.45\textwidth}
      \includegraphics[width=\textwidth]{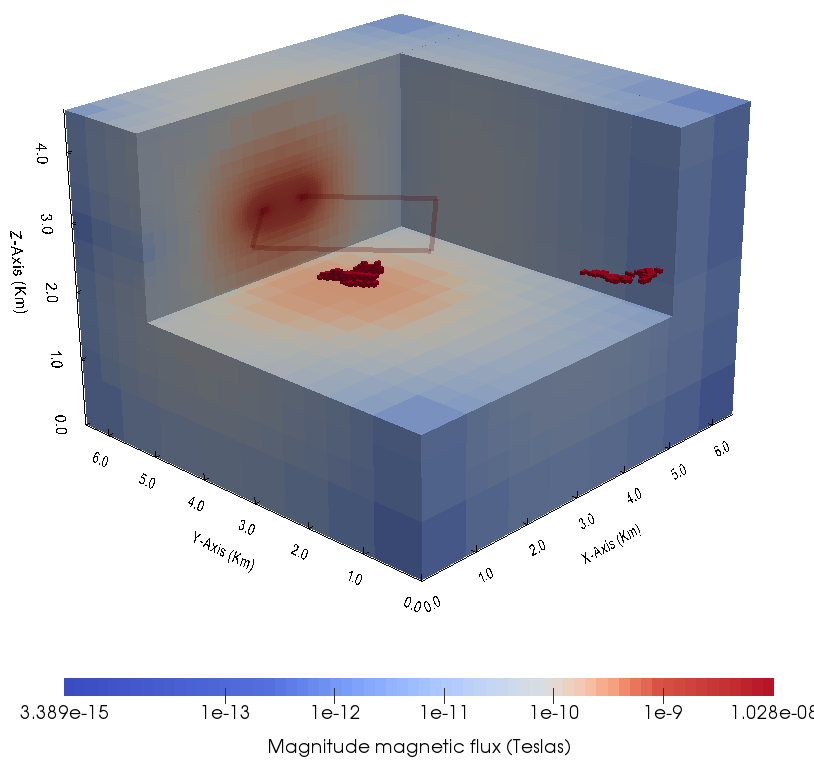}
    \end{subfigure}\\   
    \bottomrule 
  \end{tabular}    
  	\caption{\small{Magnitude of the magnetic flux for the large-loop EM survey setting.  First and second rows show the simulation results for the frequency of 1 and 20 Hz, respectively. The left-hand panel shows the reference solution obtained from solving the reference model on a finer OcTree mesh.  The right-hand panel shows the solution obtained using the proposed upscaling framework with 8 padding cells as the size of every extended local domain. }}\label{fig:simRes}
\end{figure}

The coarse conductivity models constructed using the upscaling framework were used to simulate magnetic data for the large-loop EM survey.  To do so, we used the MFV discretization method on the described coarse OcTree mesh \cite{Haber2007a,Horesh2011}.  This discretization yields linear systems of equations with 169,892 edge degrees of freedom (DOF).  The systems were solved using the direct solver MUMPS \cite{Amestoy2001}. The magnitude of the magnetic flux, denoted as $|\B|$, obtained for the frequencies of 1 and 20~Hz are shown in the right-hand panel of Figure \ref{fig:simRes}. 

To assess the accuracy of the upscaling framework, our simulation results were compared to the simulation results obtained by solving the conductivity model on the finer OcTree mesh (Figure \ref{fig:referenceModel}).  Once again, we used MFV for the discretization procedure, which yields systems with roughly 1.5 millions DOF, and solved the systems using MUMPS.  The magnitude of the resulting magnetic flux, $|\B|$, for the frequencies of 1 and 20~Hz is shown in the left-hand column of Figure \ref{fig:simRes}.

Note that the linear system for the fine-mesh problem has roughly 10 times more DOF than the coarse-mesh system.  Using a direct method to solve  the fine-mesh system, which scales as a cubic power of the variables, translates to a computation time that is roughly 30 times longer than that of the coarse-mesh system (cf. \cite{Golub1996}).

To evaluate the quality of the magnetic data produced by the coarse conductivity models constructed with the proposed upscaling framework, we also carried out simulations using additional coarse-mesh conductivity models that were constructed using average-based upscaling approaches.   The comparison among the various magnetic data was done as follows. 
We first constructed each additional coarse model by using arithmetic, geometric and harmonic averaging of the fine-mesh conductivity inside each coarse cell of the coarse OcTree mesh shown in Figure \ref{fig:CoarseLalorSetting}.  We then simulated magnetic data for the large-loop EM survey setting for each additional coarse model and individual frequency. 

To compute the error between the different upscaling techniques, we use secondary fields.
The secondary field was computed as
\begin{equation}
	\Delta \B(\Sigma) =  \B(\Sigma) - \B(\sigma_{\sf b}),
\end{equation}
where the background surface conductivity, $\sigma_{\sf b}$, is 0.01 S/m, and $\Sigma$ denotes an electrical conductivity model.

The relative error is given by
\begin{equation}
 \delta_{\Delta \B}  = \frac{|| ~|\Delta \B^h| - |\Delta \B^H|~ ||_\infty}{||\Delta \B^H||_\infty }\times 100,
\end{equation}
where the superscripts $h$ and $H$ denote dependence of the fine and coarse OcTree mesh used to compute $\Delta \B$, respectively.

\begin{table}[!htb]
\centering
\caption{Relative error for the secondary magnetic flux data in magnitude, $\delta_{\Delta \B}$, resulting from forward modeling using five coarse-mesh conductivity models and two frequencies.  \label{tbl:relErrorsBs}}
\begin{subtable}{\textwidth}
\centering
\begin{tabular}{  l c c c}
	\toprule
	\multicolumn{1}{p{4.0 cm}}{\centering Coarse-mesh \\ conductivity model } & 
	\multicolumn{1}{p{4.4 cm}}{\centering $\delta_{\Delta \B}$ at 1 Hz \\ (percent)}   & 
	\multicolumn{1}{p{4.5 cm}}{\centering $\delta_{\Delta \B}$ at 20 Hz \\ (percent)}  \\
	\midrule
 	Arithmetic 		  			 &  9.028  & 1.534  \\
  	Geometric 		  			 &  8.990  & 0.594  \\ 
  	Harmonic 		  			 &  8.987  & 0.403  \\ 
    Upscaling (4 padding cells)  &  8.989  & 0.532  \\ 
    Upscaling (8 padding cells)  &  8.991  & 0.383  \\
 	\bottomrule
 \end{tabular}
\end{subtable}
\end{table}

Examining the results, we see that while upscaling has an effect when using 
20~Hz data the effect is smaller when considering the 1~Hz data. This should not come as
a surprise: the fields at 1~Hz are mainly in the real component, and are less sensitive to fine-scale variations in conductivity than at 20~Hz. 
Indeed, for the magnetostatic case, the magnetic fields are conductivity independent. 
As a result, the error that we observe in this frequency is mostly due to discretization error.
However, when considering 20~Hz, the effect of using an appropriate averaging scheme is more evident and in fact, our averaging scheme does better than other averaging schemes.
Nonetheless, for this case, it is surprising to see how well simple harmonic averaging
did. We suspect that, unlike the 1D case where harmonic averaging performs poorly, the 
3D conductivity model requires averaging over a much smaller, local area which leads to
a much smaller difference in the data.

\section{Conclusions} \label{sec:conclusions}

In this paper, we develop a numerical upscaling framework to construct coarse-mesh electrical conductivity models based on prescribed fine-mesh ones for a broad range of quasi-static EM geophysical problems in frequency domain.  In practice, simulating these types of problems is computationally expensive; they often consider highly heterogeneous geologic media that require a very large and fine mesh to be discretized accurately. In the proposed framework, we pose upscaling as a parameter estimation problem.  Thus, a coarse-mesh conductivity model is obtained by solving an optimization problem for each coarse-mesh cell (possibly in parallel).  The optimization criterion (i.e., upscaling criterion) can be customized to construct isotropic or fully anisotropic real or even complex upscaled quantities to approximate any of the EM responses (fields and/or fluxes) depending on the geophysical experiment of interest. In particular, different experiments use different upscaling criteria that result in different upscaled quantities. As a consequence, the framework is able to upscale arbitrary electrical conductivity structures in an effective and accurate manner.  Our 1D and 3D experiments show that the coarse-mesh models constructed with the proposed upscaling framework yield accurate approximations to the EM responses that are comparable to those obtained by using a fine mesh in the forward modeling process, and that the size of the problem can be reduced significantly. For the examples presented, the size of the coarse-mesh system solved was roughly 10\% of the fine-mesh system size, while the relative errors (in the secondary fields) were less than 5\%.  That is, the coarse conductivity models are able to emulate the behavior of the heterogeneity present in the prescribed fine-mesh conductivity model.

The upscaling framework has some disadvantages. First, it constructs an upscaled conductivity that depends significantly on the set of boundary conditions imposed to compute the synthetic data used in the upscaling criterion.  We use the set of standard bilinear decaying functions on a coarse cell as the set of boundary conditions based on the arguments described in Section 4.  We recognize that such a set of boundary conditions may not be the most appropriate for constructing accurate coarse models for all cases.  However, for the experiments presented, these boundary conditions give reasonable estimates.  Second, our method is more expensive than simple average-based upscaling methods, as we solve a local optimization problem in each coarse cell.  However, since each local problem is formulated independently of the rest (even if we extend the local domain), we can reduce the cost by using a parallel implementation of the method.  

The upscaling framework demonstrates that the construction of upscaled quantities should be specific to, and highly dependent on, the purpose of the simulation. The choices of the EM responses of interest, boundary conditions, and type of the upscaled quantity employed (i.e., isotropic or anisotropic) all influence the nature of the resulting upscaled conductivity model. As a result, for a given fine-scale conductivity structure, there is no unique upscaled model which completely describes it. The inherent non-uniqueness of the proposed upscaling process is a feature that can provide insights into the behavior of the EM responses due to the presence of heterogeneous geologic materials.

\section{Acknowledgments}
The authors would like to thank Douglas Oldenburg for his insightful comments and discussions on the upscaling method, Dominique Fournier for generously providing the data that we used to construct the synthetic conductivity model used for the 3D numerical experiments presented, and Gudni Rosenkjar for his assistance with the software Paraview that we use to visualize the 3D numerical results. 
We also thank the anonymous referees for their thoughtful and constructive comments that help to improve this paper.
The funding for this work is provided through the University of British Columbia's Four-Year-Fellowship program and the Vanier Canada Graduate Scholarships Program.
\clearpage
\newpage
\bibliographystyle{plain}
\bibliography{library,biblio,lindseyRefs}

\begin{thebibliography}{10}

\bibitem{Amestoy2001}
P.~Amestoy, I.~Duff, J.~Y. L'Excellent, and J.~Koster.
\newblock {MUMPS: a general purpose distributed memory sparse solver}.
\newblock {\em Int. Work. Appl. Parallel Comput. Springer Berlin Heidelb.},
  pages 121--130, 2001.

\bibitem{Avdeev2005}
D.~Avdeev.
\newblock {Three-dimensional electromagnetic modelling and inversion from
  theory to application}.
\newblock {\em Surv. Geophys.}, 26(6):767--799, 2005.

\bibitem{BerrymanHoversten2013}
J.~G. Berryman and G.~M. Hoversten.
\newblock {Modelling electrical conductivity for earth media with macroscopic
  fluid-filled fractures}.
\newblock {\em Geophys. Prospect.}, 61(2):471--493, mar 2013.

\bibitem{Borner2009}
R.~U. B{\"{o}}rner.
\newblock {Numerical modelling in geo-electromagnetics: advances and
  challenges}.
\newblock {\em Surv. Geophys.}, 31(2):225--245, oct 2010.

\bibitem{caudillo2014}
L.~A. Caudillo-Mata, E.~Haber, L.~J. Heagy, and D.~W. Oldenburg.
\newblock {\em {Numerical upscaling of electrical conductivity: A problem
  specific approach to generate coarse-scale models}}, chapter 130, pages
  680--684.

\bibitem{Durlofsky1998}
L.~J. Durlofsky.
\newblock {Coarse scale models of two phase flow in heterogeneous reservoirs:
  volume averaged equations and their relationship to existing upscaling
  techniques}.
\newblock {\em Comput. Geosci.}, 2(2):73--92, mar 1998.

\bibitem{Durlofsky2003}
L.~J. Durlofsky.
\newblock {Upscaling of geocellular models for reservoir flow simulation: a
  review of recent progress}.
\newblock {\em 7th Int. Forum Reserv. Simul. B{\"{u}}hl/Baden-Baden, Ger.},
  pages 1--58, 2003.

\bibitem{Efendiev2009}
Y.~Efendiev and T.~Y. Hou.
\newblock {\em {Multiscale finite element methods: theory and applications}}.
\newblock Springer New York, 2009.

\bibitem{Farmer2002}
C.~L. Farmer.
\newblock {Upscaling: a review}.
\newblock {\em Int. J. Numer. Methods Fluids}, 40(1-2):63--78, 2002.

\bibitem{Farquharson2003}
C.~Farquharson, D.~Oldenburg, and P.~S. Routh.
\newblock {Simultaneous 1D inversion of loop-loop electromagnetic data for
  magnetic susceptibility and electrical conductivity}.
\newblock {\em Geophysics}, 68(6):1857--1869, 2003.

\bibitem{Golub1996}
G.~H. Golub and C.~F.~V. Loan.
\newblock {\em {Matrix Computations (Johns Hopkins Studies in Mathematical
  Sciences)(3rd Edition)}}, volume 208-209.
\newblock 1996.

\bibitem{haberBook2014}
E.~Haber.
\newblock {\em Computational Methods in Geophysical Electromagnetics}.
\newblock SIAM, Philadelphia, 2014.

\bibitem{Haber2014}
E.~Haber.
\newblock {\em {Computational Methods in Geophysical Electromagnetics}}.
\newblock Society for Industrial and Applied Mathematics, 2014.

\bibitem{Haber2007a}
E.~Haber and S.~Heldmann.
\newblock {An octree multigrid method for quasi-static Maxwell's equations with
  highly discontinuous coefficients}.
\newblock {\em J. Comput. Phys.}, 223(2):783--796, may 2007.

\bibitem{Haber2014c}
E~Haber and L~Ruthotto.
\newblock {A multiscale finite volume method for Maxwell's equations at low
  frequencies}.
\newblock {\em Geophys. J. Int.}, 199(2):1268--1277, 2014.

\bibitem{Horesh2011}
L.~Horesh and E.~Haber.
\newblock {A Second Order Discretization of Maxwell's Equations in the
  Quasi-Static Regime on OcTree Grids}.
\newblock {\em SIAM J. Sci. Comput.}, 33(5):2805--2819, 2011.

\bibitem{Hyman1998}
J.~M. Hyman and M.~Shashkov.
\newblock {Mimetic discretizations for Maxwell equations and the equations of
  magnetic diffusion}.
\newblock Technical report, United States. Department of Energy. Office of
  Energy Research, 1998.

\bibitem{Hyman1999}
J.~M. Hyman and M.~Shashkov.
\newblock {Mimetic discretizations for Maxwell's equations}.
\newblock {\em J. Comput. Phys.}, 151(2):881--909, may 1999.

\bibitem{Hyman1999a}
J.~M. Hyman and M.~Shashkov.
\newblock {The orthogonal decomposition theorems for mimetic finite difference
  methods}.
\newblock {\em SIAM J. Numer. Anal.}, 36(3):788--818, jan 1999.

\bibitem{Jin2002}
J.~Jin.
\newblock {\em {The Finite Element Method in Electromagnetics}}.
\newblock Wiley, New York, 2nd edition, 2002.

\bibitem{Kelley1999}
C.~T. Kelley.
\newblock {\em {Iterative methods for optimization}}.
\newblock Society for Industrial and Applied Mathematics, 1999.

\bibitem{Key2011}
K.~Key and J.~Ovall.
\newblock {A parallel goal-oriented adaptive finite element method for 2.5-D
  electromagnetic modelling}.
\newblock {\em Geophys. J. Int.}, 186:137--154, 2011.

\bibitem{KristenssonWellander2003}
G.~Kristensson and N.~Wellander.
\newblock {Homogenization of the Maxwell equations at fixed frequency}.
\newblock {\em SIAM J. Appl. Math.}, 64(1):170--195, 2003.

\bibitem{Lin1999}
C.~J. Lin and J.~J. Mor{\'{e}}.
\newblock {Newton's method for large bound-constrained optimization problems}.
\newblock {\em SIAM J. Optim.}, 9(4):1100--1127, jan 1999.

\bibitem{Lipnikov2004}
K.~Lipnikov, J.~Morel, and M.~Shashkov.
\newblock {Mimetic finite difference methods for diffusion equations on
  non-orthogonal non-conformal meshes}.
\newblock {\em J. Comput. Phys.}, 199:589--597, 2004.

\bibitem{MacLachlan2004}
S.~MacLachlan.
\newblock {\em {Improving Robustness in Multiscale Methods}}.
\newblock PhD thesis, University of Colorado, 2004.

\bibitem{MacLachlan2006}
S.~MacLachlan and J.~D. Moulton.
\newblock {Multilevel upscaling through variational coarsening}.
\newblock {\em Water Resour. Res.}, 42(2):W02418:1--9, 2006.

\bibitem{milton2002theory}
Graeme~W Milton.
\newblock {\em The theory of composites}, volume~6.
\newblock Cambridge University Press, 2002.

\bibitem{Monk2003}
P.~Monk.
\newblock {\em {Finite element methods for Maxwell's equations}}.
\newblock Clarendon, 2003.

\bibitem{Moulton1998}
J.~D. Moulton, J.~E. Dendy, and J.~M. Hyman.
\newblock {The black box multigrid numerical homogenization algorithm}.
\newblock {\em J. Comput. Phys.}, 108(142):80--108, 1998.

\bibitem{Nocedal2006}
J.~Nocedal and S.~J. Wright.
\newblock {\em {Numerical Optimization}}.
\newblock Springer New York, 2nd edition, 2006.

\bibitem{Oldenburg2007}
D.~W. Oldenburg and D.~A. Pratt.
\newblock {Geophysical inversion for mineral exploration: A decade of progress
  in theory and practice}.
\newblock {\em Proc. Explor.}, 2007.

\bibitem{Schwarzbach2009}
C.~Schwarzbach.
\newblock {\em {Stability of finite element discretization of Maxwell's
  equations for geophysical applications}}.
\newblock PhD thesis, University of Frieberg, Germany, 2009.

\bibitem{Shafiro2000}
B.~Shafiro and M.~Kachanov.
\newblock {Anisotropic effective conductivity of materials with nonrandomly
  oriented inclusions of diverse ellipsoidal shapes}.
\newblock {\em J. Appl. Phys.}, 87(12):8561--8569, 2000.

\bibitem{Smith2013}
R.~Smith.
\newblock {Electromagnetic Induction Methods in Mining Geophysics from 2008 to
  2012}.
\newblock {\em Surv. Geophys.}, 35(1):123--156, apr 2013.

\bibitem{Telford1990}
W.~M. Telford, L.~P. Geldart, and R.~E. Sheriff.
\newblock {\em {Applied Geophysics}}, volume 2nd.
\newblock Cambridge University Press, 2nd edition, 1990.

\bibitem{Torquato2002}
S.~Torquato.
\newblock {\em {Random heterogeneous materials: microstructure and macroscopic
  properties}}.
\newblock Springer, 2002.

\bibitem{Wang2014}
Z.~Wang, E.~G. Lim, Y.~Tang, and M.~Leach.
\newblock {Medical Applications of Microwave Imaging}.
\newblock {\em Sci. World J.}, 2014(Article ID 147016):1--7, 2014.

\bibitem{Ward1988}
S.~H. Ward and G.~W. Hohmann.
\newblock {Electromagnetic theory for geophysical applications}.
\newblock In Misac~N. Nabighian, editor, {\em Electromagn. methods Appl.
  Geophys.}, volume~1, chapter~4, pages 130--311. Society of Exploration
  Geophysicists, 1988.

\bibitem{Weaver1999}
J.~T. Weaver.
\newblock {Numerical Modelling in Electromagnetic Induction}.
\newblock In K.~K. Roy, S.~K. Verma, and K.~Mallick, editors, {\em Deep
  Electromagn. Explor.}, chapter~19, pages 299--363. Springer, 1999.

\bibitem{Wen1996}
X.~H. Wen and J.~J. G{\'{o}}mez-Hern{\'{a}}ndez.
\newblock {Upscaling hydraulic conductivities in heterogeneous media: An
  overview}.
\newblock {\em J. Hydrol.}, 183(1):ix -- xxxii, 1996.

\bibitem{wynne1994athabasca}
Desmond~A. Wynne, Michelle Attalla, Tim Berezniuk, Habtemicael Brulotte,
  Darrell~K. Cotterill, Rudy Strobl, and Daryl~M. Wightman.
\newblock {Athabasca Oil Sands Database}.
\newblock Technical report, Alberta Geological Survey, Edmonton,AB, 1994.

\bibitem{Yang2014}
D.~Yang, D.~Fournier, and D.~W. Oldenburg.
\newblock {3D Inversion of EM data at Lalor mine: in pursuit of a unified
  electrical conductivity model}.
\newblock In {\em Explor. Deep VMS Ore Bodies Hudbay Lalor Case Study}, pages
  1--4. BC Geophysical Society, 2014.

\bibitem{Zhdanov2010}
M.~S. Zhdanov.
\newblock {Electromagnetic geophysics: Notes from the past and the road ahead}.
\newblock {\em Geophysics}, 75(5):75A49--75A66, sep 2010.

\end{thebibliography}
\end{document}